\documentclass[letterpaper, 10 pt, conference]{ieeeconf}  

\IEEEoverridecommandlockouts                              

\usepackage{xr-hyper}
\usepackage{upgreek}
\usepackage{color}
\usepackage{colortbl}
\usepackage{empheq}
\usepackage{enumerate}
\usepackage{paralist}

\usepackage{times}
\usepackage{amsmath,amsfonts,bm}
\usepackage[ruled,algo2e,vlined]{algorithm2e}
\usepackage{algpseudocode} 
\usepackage{algorithm}
\usepackage{graphicx}
\usepackage{subcaption}
\usepackage{comment}
\usepackage{xspace}
\usepackage{booktabs}
\usepackage{caption} 
\usepackage{hyperref}
\usepackage{multirow} 
\usepackage{makecell}

\newcommand{\bu}{\bm u}
\newcommand{\bR}{\mathbb {R}}
\newcommand{\bx}{\bm x}

\newcommand{\rmd}{\,\mathrm{d}}
\newcommand{\ie}{\textit{i.e.,}}
\newcommand{\eg}{\textit{e.g.,}}

\newcommand{\LQR}{\text{LQR}}
\newcommand{\NN}{\text{NN}}
\newcommand{\QR}{\text{QR}}

\newcommand{\tmin}{\text{min}}
\newcommand{\tmax}{\text{max}}

\DeclareMathOperator*{\argmin}{arg\,min}

\newcommand{\revise}[1]{#1}

\newcommand{\tabincell}[2]{\renewcommand{\arraystretch}{0.65}\begin{tabular}{@{}#1 @{}}\addlinespace[0.8ex]#2\end{tabular}\renewcommand{\arraystretch}{1}}

\newcommand{\DAgger}{\textsc{DAgger}}

\title{\LARGE \bf
Learning Free Terminal Time Optimal Closed-loop Control of Manipulators
}

\PassOptionsToPackage{hyperfootnotes=false}{hyperref}

\author{Wei Hu$^{1,2,*}$, Yue Zhao$^{3,*}$, Weinan E$^{3}$, Jiequn Han$^{4}$, Jihao Long$^{2,\dag}$
\thanks{*These authors contributed equally to this work.}
\thanks{\dag Corresponding author.  {\tt\small longjh1998@gmail.com}}
\thanks{$^{1}$AI for Science Institute, Beijing, China}
\thanks{$^{2}$Institute for Advanced Algorithms Research, Shanghai, China}
\thanks{$^{3}$Center for Data Science, Peking University, Beijing, China}%
\thanks{$^{4}$Flatiron Institute, New York, USA}%
}

\begin{document}

\maketitle
\thispagestyle{empty}
\pagestyle{empty}

\begin{abstract}%
This paper presents a novel approach to learning free terminal time closed-loop control for robotic manipulation tasks, enabling dynamic adjustment of task duration and control inputs to enhance performance. We extend the supervised learning approach, namely solving selected optimal open-loop problems and utilizing them as training data for a policy network, to the free terminal time scenario. Three main challenges are addressed in this extension. First, we introduce a marching scheme that enhances the solution quality and increases the success rate of the open-loop solver by gradually refining time discretization. Second, we extend the QRnet in \cite{QRnet} to the free terminal time setting to address discontinuity and improve stability at the terminal state. Third, we present a more automated version of the initial value problem (IVP) enhanced sampling method from previous work \cite{IVP_Enhanced} to adaptively update the training dataset, significantly improving its quality. By integrating these techniques, we develop a closed-loop policy that operates effectively over a broad domain with varying optimal time durations, achieving near globally optimal total costs. The videos are available at \textcolor{blue}{\href{https://deepoptimalcontrol.github.io/FreeTimeManipulator}{https://deepoptimalcontrol.github.io/FreeTimeManipulator}}.
\end{abstract}

\section{Introduction}
Classical manipulator control strategies, such as proportional–integral–derivative~(PID), sliding mode, and Lyapunov control, are primarily engineered for tracking accuracy and robustness, often achieved by leveraging properties like passivity and feedback linearization of the manipulator
\cite{historial-control-manipulator2022}.
These methods track trajectories from high-level path planners, which often ignore the manipulators’ intrinsic dynamics, resulting in suboptimal control outcomes.
Model predictive control (MPC) offers an alternative with the potential for better performance, but it still compromises long-term optimality and requires real-time online optimization \cite{MPC1,MPC2}.
In recent years, deep reinforcement learning (RL) has gained popularity in robotic manipulation, primarily functioning as a high-level path planner that generates position or velocity commands \cite{rl-review, rl-visual-pcontrol,rl-vcontrol}. However, state-of-the-art RL algorithms still fall short in achieving true optimality as the low-level control \cite{hu2022solving,zhao2024offline}.
To break these limitations, this paper presents a methodology that trains a neural network controller using a free terminal time optimal control formulation, integrating runtime optimization with other practical objectives. This approach simultaneously and correctly adjusts the remaining task time and low-level control inputs to enhance performance. 

The concept of free terminal time control has proven valuable not only in robotics but also in fields like space trajectory optimization~\cite{liu2019fuel} and chemical processing~\cite{palanki1994optimal}. However, solving the numerical challenges of the free terminal time open-loop optimal control problem remains difficult. Specifically, this problem requires optimizing the time duration along with high-dimensional states and controls, all within the constraints of highly non-linear and non-convex dynamics. The complexity further increases with closed-loop control strategies that must adapt to various initial states. Traditional mesh-based discretization methods for solving the Hamilton-Jacobi-Bellman (HJB) equation~\cite{bardi1997optimal,cristiani2010initialization,falcone2014semi} suffer from the curse of dimensionality. Moreover, the optimal value or control function often exhibits discontinuities or non-differentiability at the terminal state~\cite{bardi1997optimal}.

Recently, deep neural networks have become prominent in approximating high-dimensional optimal closed-loop control policies, addressing the curse of dimensionality problem~\cite{han2016,silver2017mastering,Adaptive_BVP_HJB}. For an in-depth discussion on neural network-based methods in learning optimal closed-loop control, see \cite{zhao2024offline}. While most research has focused on fixed or infinite time horizon controls, free terminal time controls have been less explored. An exception is \cite{freetimeac2023}, which introduced an iterative method using neural networks to approximate value and policy functions for solving the free terminal time HJB equation. Although effective in low-dimensional scenarios, a gap remains between the learned and optimal policies in the 7-dimensional attitude control of spacecraft.

In this paper, we extend the supervised learning approach to learn a free terminal time optimal controller for robotic manipulators. Instead of solving the HJB equation, we generate optimal open-loop trajectories for various initial states using numerical solvers to create a dataset. We then train a policy network to fit this dataset. To this end, we overcome several algorithmic challenges uniquely associated with free terminal time optimal control problems. In the open-loop solver, we introduce a marching scheme that progresses from coarse to fine-time discretization. This approach notably enhances the solution quality and increases the convergence rate of the solver.
On the supervised learning side, we address two key obstacles. First, network-based control often suffers from stability issues at the terminal state, which are further complicated by the discontinuity in the optimal closed-loop control at the terminal state in our problem. To overcome these challenges, we extend the QRnet proposed in \cite{QRnet} to the free terminal time setting to capture this discontinuity and significantly improve stability. Second, the dataset quality greatly affects the final learned policy since the network can perform poorly on unrepresented states. Prior work~\cite{IVP_Enhanced} proposed an initial value problem (IVP) based enhanced sampling method to systematically update the dataset for fixed terminal time problems. However, that procedure requires pre-specified time grids to select samples. We propose an adaptive and general IVP sampling method for enhanced sampling of free terminal time problems,
which automatically identifies time grids on a per-trajectory basis.
By combining these techniques, we obtain a closed-loop policy working over a large domain, where the optimal time duration varies substantially, and our policy achieves near globally optimal total cost.

\section{Preliminaries}\label{sec:prob}

\subsection{Optimal Control Problem}
A deterministic controlled dynamical system in general can be written as
\begin{equation}\label{eq:dyanmic_system}
  \dot{\bm{x}}(t) = \bm f( \bm{x}(t), \bm{u}(t)), \; t\in[0,t_f],\; 
     \bm{x}(0) = \bm{x}_0,\; \bm{x}(t_f) = \bm{x}_f,
\end{equation}
where $\bm{x}(t) \in \bR^n$ represents the state, $\bm{u}(t) \in \mathcal{U} \subset \bR^m$ denotes the control with $\mathcal{U}$ being the set of admissible controls, $t_f \ge 0$ is the terminal time, $\bm{f}\colon\bR^n\times \mathcal{U} \rightarrow \bR^n$ is a smooth function describing the dynamics, and $\bm{x}_0, \bm{x}_f \in \bR^n$ denote the initial and terminal states. In this work, we assume $\bx_f$ is fixed.
Given a $t_f$, 
we call a path $\{\bm{x},\bm{u}\}_{t=0}^{t_f}$ feasible if it satisfies \eqref{eq:dyanmic_system}, and $\mathcal{P}(t_f,\bm{x}_0)$ denotes the set of all such paths.
The total cost is:
\begin{equation}\label{eq:cost}
  \mathcal{C}[\bm{x},\bm{u},t_f] = \int^{t_f}_{0} L(\bx, \bu)\rmd t,
\end{equation}
where $L:\bR^n\times \mathcal{U}\rightarrow \bR$ is the running cost.

The \emph{free terminal time open-loop problem}, given initial state $\bx_0$,  is to find solution of
\begin{equation}
\label{eq:free_open_problem}
\underset{t_f,\{\bm{x},\bm{u}\}\in\mathcal{P}(t_f,\bm{x}_0)}{\operatorname{minimize}}
 \mathcal{C}[\bm{x},\bm{u},t_f]=\colon V(\bm{x}_0).
\end{equation}
We use $t_f^*(\bm{x}_0)$, $\bm{x}^*(t;\bm{x}_0)$, and $\bm{u}^*(t;\bm{x}_0)$ to denote the optimal terminal time, state, and open-loop control associated with the initial state $\bm{x}_0$, respectively, highlighting their dependencies on the initial state. 
Since the terminal state $\bx_f$ will not vary throughout the paper, we omit the dependence on it.
For simplicity, we assume \eqref{eq:free_open_problem} is well-posed, \ie~the solution always exists and is unique. \revise{However, we acknowledge that in some cases, solutions may not be unique. In such scenarios, alternative approaches, such as probabilistic methods or diffusion-based policies~\cite{chi2023diffusion}, may be more appropriate. Exploring these directions is an interesting avenue for future work.}

In contrast to the open-loop control being a function of time, the closed-loop control $\bm{u}\colon\mathbb{R}^n \rightarrow \mathcal{U}$ of the free terminal time problem is a function of the state $\bm x$. A closed-loop control  induces a family of the open-loop controls with different $\bm{x}_0$, by solving the \emph{initial value problem} (IVP) denoted as $\text{IVP}(\bm x_0, \bm{u})$:
\begin{equation}\label{eq:IVP}
\begin{dcases}
\dot{\bm{x}}^{\bm u}(t;\bm{x}_0) = \bm{f}(\bm{x}^{\bm u}(t;\bm{x}_0),\bm{u}(\bm{x}^{\bm u}(t;\bm{x}_0))),\,
\bm{x}^{\bm u}(0;\bm{x}_0) = \bm{x}_0, \\
t_f^{\bm{u}}(\bm{x}_0) = \inf\{t \ge0\colon\bm{x}^{\bm{u}}(t,\bm{x}_0) = \bm{x}_f\} \,(\inf \emptyset = + \infty),\\
\bm{u}(t;\bm{x}_0) = \bm{u}(\bm{x}^{\bm u}(t;\bm{x}_0)).
\end{dcases}
\end{equation}
For simplicity, we use the same symbol to represent both the closed-loop control function and the induced open-loop controls, discernible from the context.
Classical optimal control theory (see, \eg \,\cite{bardi1997optimal}) asserts:  
\begin{inparaenum}
    \item [1)]there exists  a closed-loop optimal control function $\bm{u}^*\colon\mathbb{R}^n \rightarrow \mathcal{U}$ such that for any initial state $\bm{x}_0$:
$
\bm{u}^*(t;\bm{x}_0) = \bm{u}^*(\bm{x}^{\bm u^*}(t;\bm{x}_0)),
$
implying that the family of open-loop optimal controls with all possible initial states can be induced from the closed-loop optimal control function; 
    \item [2)]the value function in \eqref{eq:free_open_problem} is the viscosity solution of the HJB equation 
    \begin{equation}
    \label{eq:HJB}
        \min_{\bu \in \mathcal{U}}\{L(\bm{x},\bm{u}) + \nabla V(\bm{x}) \cdot \bm{f}(\bm{x},\bm{u})\} = 0,\; V(\bm{x}_f) = 0,
    \end{equation}
    and, when $\nabla V(\bx)$ exists,
    $
        \bm u^*(\bm{x}) = \argmin_{\bm u \in \mathcal{U}}\{L(\bm{x},\bm{u}) + \nabla V(\bm{x})\cdot \bm{f}(\bm{x},\bm{u})\}
    $.
\end{inparaenum}

\subsection{Manipulator}\label{Sec: Manipulator}
As a prototypical example, we consider a 7-DoF manipulator with torque control, 
the KUKA LWR iiwa 820 14 \cite{kuka-platform}, whose dynamics can be written as
$
    \dot{\bm x} = \bm f(\bm x, \bm u) = (\bm v, \bm a(\bm x, \bm u))
$, where the state $\bm x = (\bm q, \bm v) \in \bR^{14}$ is described by the joint angles $\bm q \in \bR^7$ and velocities $\bm v = \dot{\bm q} \in \bR^7$. The control inputs, $\bm u \in \bR^7$, denote the applied torques at the joints.
The joint accelerations, $\ddot{\bm q} = \bm a(\bm x, \bm u) \in \bR^7$, are governed by the manipulator dynamics \cite{MR2676219}:
\begin{equation}
    M(\bm q)\bm a +  C(\bm q, \dot{\bm q})\dot{\bm q} + \bm g(\bm q)= \bm u,
\end{equation}
where $M(\bm q)$ is the generalized inertia matrix \revise{which is positive definite~\cite{MR1300410}, and} $C(\bm q,\dot{\bm q})\dot{\bm q}$ represents the \emph{centrifugal} forces and \emph{Coriolis} forces, and $\bm g(\bm q)$ is the generalized gravity.
Problem \eqref{eq:free_open_problem} defines the optimal reaching problem of the manipulator. The desired terminal state $\bx_f$ is $(\bm q_f, 0)$, a state of zero velocity. 
With $\bu_f=\bm g(\bm{q}_f)$,
the state-control pair $(\bx_f, \bu_f)$ forms a static point where $\bm f(\bx_f, \bu_f) = 0$. The running cost is 
\begin{equation}
    \label{eq:running_cost_used}
    L(\bx, \bu)=r_t+r_u\|\bu - \bu_f\|^2 + r_a\|\bm a(\bx, \bu)\|^2
\end{equation}
where scalars $r_t$, $r_u$, and $r_a$ are weights for time, control, and acceleration, respectively.

\section{Learning Closed-loop Controller from Open-loop Solutions}

At a high level, our method is an instance of exploration-labeling-training (ELT) algorithms \cite{zhang2018reinforced, han2021machine}, extended to free terminal time problems. The supervised learning approach entails generating a dataset $\mathcal{D}$ of open-loop solutions (\emph{labeling}) and learning a closed-loop controller $\hat{\bm{u}}$ to approximate these solutions (\emph{training}). Additionally, we enhance the IVP sampling method \cite{IVP_Enhanced} for \emph{exploration}.

Building upon this framework, we introduce enhancements in three key aspects.
First, Section~\ref{sec:solver} presents an iterative solver for free terminal time problems, which efficiently and robustly generates desired open-loop solutions, thus creating datasets of data triples:
\begin{equation}
    \label{eq:dataset}
    \mathcal{D} = \{
    (\bx^{(i)}, \bu^{(i)}, t_f^{(i)}): \bu^{(i)} = \bu^*(\bx^{(i)}), t_f^{(i)} = t_f^*(\bx^{(i)})\}.
\end{equation}
The dataset $\mathcal{D}$ will be used to train a policy network $\hat{\bu}$ with the supervision loss:
\begin{equation}
\label{eq:u_loss}
\underset{\theta}{\operatorname{minimize}}\quad \frac{1}{|\mathcal{D}|} \sum_{i=1}^{|\mathcal{D}|}\left\|\hat{\bm{u}}(\bm{x}^{(i)} ; \theta) - \bm{u}^{(i)}\right \|^2.
\end{equation}
Second, we adapt QRnet~\cite{QRnet} to handle the discontinuity at the terminal state in free terminal time scenarios in Section~\ref{sec:qrnet}. Finally, we develop the IVP enhanced sampling technique~\cite{IVP_Enhanced} with automatic resampling times in Section~\ref{sec:IVP}.

\subsection{Iterative Solver for Free Terminal Time Data Generation}\label{sec:solver}
For generating optimal open-loop solutions, we can rewrite the problem \eqref{eq:free_open_problem} as follows:
\begin{equation}
\label{eq:free_open_problem2}
 \underset{t_f}{\operatorname{minimize}} \,\,V[\bx_0;t_f] \coloneqq  \underset{t_f}{\operatorname{minimize}}\underset{\{\bm{x},\bm{u}\}\in\mathcal{P}(t_f, \bx_0)}{\operatorname{minimize}}\mathcal{C}[\bm{x},\bm{u},t_f].
\end{equation}
Here, the inner optimization over $\mathcal{P}(t_f, \bx_0)$ defines the \emph{fixed terminal time open-loop problem}.
The solution to the fixed terminal time problem is denoted by $\bu^*(t;\bx_0, t_f)$ for the control and $\bx^*(t;\bx_0,t_f)$ for the state.
With the optimal open-loop solution $\bx^*(t;\bx_0, t_f^*(\bx_0))$, the optimal terminal times along the trajectory are $t_f^*(\bx^*(t;\bx_0, t_f^*(\bx_0)) = t_f^*(\bx_0) - t$ for any $t\in[0, t_f^*(\bx_0)]$.

For numerical considerations, we adapt the terminal constraints $\bx(t_f)=\bx_f$ into a terminal cost function, defined as $M(\bx) = r_f|\bx-\bx_f|^2$. Consequently, the cost function $\mathcal{C}[\bx, \bu, t_f]$ is expressed as $\int_0^{t_f} L(\bx,\bu)\rmd t + M(\bx)$. Choosing a relatively high weighting factor $r_f$ ensures that the state at terminal time $\bx(t_f)$ approaches the desired terminal state $\bx_f$ with adequate accuracy.

To address the free terminal time open-loop problem, we propose an iterative solver.
Starting with an initial state $\bx_0$ and an initial estimation of the terminal time $t_f^0$, each iteration $k$ of the solver involves two main steps:\vspace{-2pt}
\begin{enumerate}
    \setlength{\itemsep}{0pt}
    \setlength{\parskip}{0pt}
    \item \textit{Open-Loop Solver:} solve the fixed terminal time open-loop problem using the terminal time \( t_f^{k-1} \) from the previous iteration, yielding the optimal control \( \bu^*(t; \bx_0, t_f^{k-1}) \) and the corresponding state trajectory \( \bx^*(t; \bx_0, t_f^{k-1}) \);
    \item \textit{Terminal Time Update:} compute the gradients \( \rmd\mathcal{C}/\rmd t_f \) at \( t_f^{k-1} \) utilizing \( \bu^*(t; \bx_0, t_f^{k-1}) \) and \( \bx^*(t; \bx_0, t_f^{k-1}) \), and then apply a gradient descent or quasi-Newton method to get updated \( t_f^k \).
\end{enumerate}
We will further elucidate these steps in the subsequent paragraphs.
Refer to Algorithm 1 in Appendix~\ref{Sec:Alg} for a comprehensive summary.

\paragraph*{A Robust Open-Loop Solver}
Differential dynamic programming (DDP, \cite{DDP}) is recognized as an efficient algorithm for solving fixed terminal time optimal control problems.
However, it faces challenges like poor convergence or convergence to suboptimal local minima in highly non-convex, nonlinear problems, especially when dealing with a vast number of unknowns.
One strategy to alleviate these challenges is to build a well-constructed initial guess.
To this end, we adopt a systematic approach to incrementally enlarge the problem size so that solutions from smaller size problems provide good initial guesses for larger size problems.
This is achieved through a sequence of increasing discretization time step counts, denoted as ${n_1, n_2, \cdots, n_J}$. The procedure at iteration $k$ is as follows:
\begin{inparaitem}
\item[1)] initially for $j=1$, the problem is solved with $n_1$ time steps using initial guess from iteration $k-1$ with $n_J$ time steps. Regarding the difference in terminal times, time is scaled as $(t_f^{k}/t_f^{k-1}) t$. If $k=1$, we adopt the zero function as the initial guess;
\item[2)] for each subsequent step $j \ge 2$, the open-loop problem is solved with $n_j$ time steps, employing the solution from the preceding $n_{j-1}$ time steps problem as the initial guess.
\end{inparaitem}

By progressively scaling the problem size, this method ensures more precise and stable solutions in complex control scenarios.
For our implementation, we employ the standard DDP algorithm as provided in the Crocoddyl library \cite{mastalli20crocoddyl}.
The numerical computation of the dynamics is conducted using the Pinocchio library \cite{carpentier2019pinocchio}.

\paragraph*{Terminal Time Update}
In our iterative approach, the terminal time $t_f$ is updated using both gradient descent and quasi-Newton methods.
For the first several iterations, we employ gradient descent updates for $t_f$, continuing until the gradient diminishes to a relatively small value.
Subsequently, we switch to the quasi-Newton method.
This strategy mitigates the limitations inherent in each method: quasi-Newton methods necessitate a more precise initial estimation and cannot differentiate between local minima and maxima. Sole reliance on gradient descent, on the other hand, may cause oscillations around a terminal time where the gradient is large.
In contrast, the quasi-Newton update significantly improves the accuracy of $t_f$ in later iterations.

\begin{figure*}[ht]
    \centering
    \includegraphics[width=1.0\textwidth]{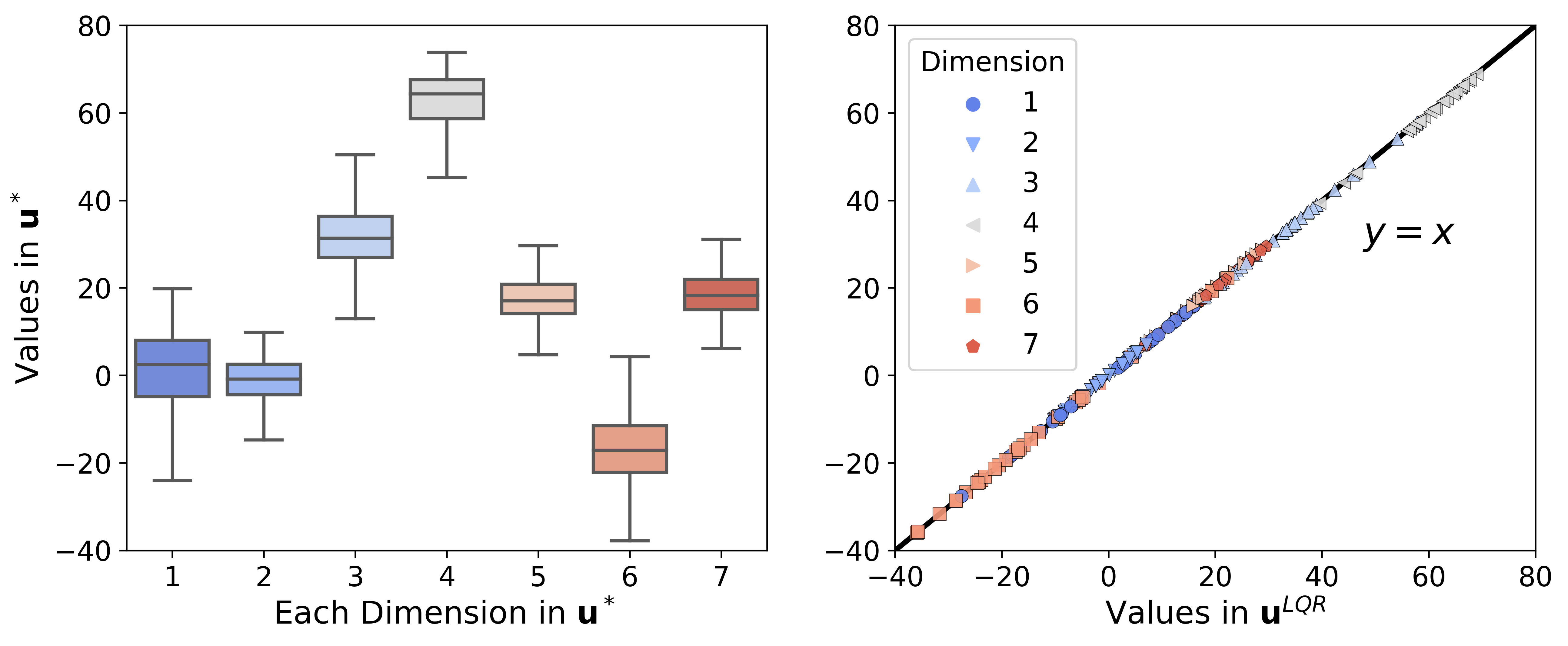}
    \caption{The left figure shows the distribution of final control from the open-loop control $\bm u^*(t_f^*(\bm x_0);\bm x_0)$ for varying initial states $\bm x_0$. Although the states $\bm x^*(t_f^*(\bm x_0);\bm x_0)$ are uniformly close to $\bx_f$ within a 0.001 Euclidean norm, there is a notable variation in optimal controls across these states in all 7 dimensions. This highlights the divergence of $\bm u^*(\cdot)$ around $\bx_f$. The right figure compares the values of $\bm{u}^{\mathrm{LQR}}$ and $\bm u^*$ at $\bm x^*(t_f^*(\bm x_0); \bm x_0)$, the final state in each open-loop solution. The average difference between them, calculated using the Euclidean norm, is 0.0768, much smaller than the norm of the original control.}
    \label{fig:terminal_control_lqr}
\end{figure*}

The gradient of the total cost $\mathcal{C}$ with respect to $t_f$ is expressed as \cite{MR446628}: 
\begin{equation}
\label{eq:graident_t}
\begin{aligned}
& \frac{\rmd}{\rmd t_f} \underset{{\bm{x},\bm{u}}\in\mathcal{P}{t_f}}{\operatorname{minimize}}\,\,\mathcal{C}[\bm{x},\bm{u},t_f]\\ 
&= \int_{0}^{t_f} H(t,\bx^*(t;\bx_0,t_f), \bu^*(t;\bx_0,t_f);t_f) \rmd t,
\end{aligned}
\end{equation}
where $H(t,\bx,\bu;t_f)=L(\bx,\bu) + \nabla V(\bx;t_f-t) \cdot \bm{f}(\bx,\bu)$.
The update scheme is detailed next.
Let $\delta t_f^{k}$ and $\rmd t_f^{k}$ represent the update direction and the gradient at iteration $k$, respectively.
The update direction for both methods is defined as:
\begin{inparaitem}
    \item[1)] the gradient descent method sets the update direction equal to the gradient, \ie~$\delta t_f^{k} = \rmd t_f^{k}$;
    \item[2)] the quasi-Newton method approximates the inverse of the Hessian matrix to refine the update direction:
    $
        \delta t_f^{k} = (t_f^{k} - t_f^{k-1})/(\rmd t_f^{k} - \rmd t_f^{k-1})\rmd t_f^{k}.
    $
\end{inparaitem}
The terminal time is then updated as follows:
$
t_f^{k+1} = t_f^{k} - \alpha^{k} \delta t_f^{k},
$
where $\alpha^k\in (0,1]$ is a step size parameter chosen to ensure that the magnitude of the update does not exceed a fraction $\alpha$ of the previous terminal time, with $\alpha \in (0,1]$ serving as a safeguard against excessively large updates.

Finally, to enhance the consistency of the data, we standardize the time step size of the solutions associated with different $\bx_0$ to a unified small value, denoted as $\Delta t$. This is achieved by rounding the converged terminal time $t_f^k$ to its nearest multiple of $\Delta t$. Subsequently, we solve the fixed terminal time open-loop problem, adjusted for the slightly altered terminal time.
The solution obtained from this process is then referred to as the numerical solution for the free terminal time problem.

\subsection{Handling Discontinuity at Terminal State}\label{sec:qrnet}
The problem in Section \ref{sec:prob} presents an additional challenge compared to the standard one: its optimal closed-loop solution is discontinuous at $\bx_f$. 
The non-smoothness of the value function $V$ can be seen by contradiction. If $V$ is differentiable within an open set around $\bx_f$, then $V(\bx_f)=0$ and $V(\bx)\geq 0$ for $\bx\ne\bx_f$ implies $\nabla V(\bx_f)=0$.
According to \eqref{eq:HJB}, we can derive $\min_{\bm u \in \mathcal{U}}\{L(\bm x_f,\bm u)\} = 0$,
which contradicts the strict positivity of $L$. Furthermore, our numerical results demonstrate the discontinuity of the optimal control function at $\bx_f$, as depicted in Figure \ref{fig:terminal_control_lqr} (left).

The discontinuity of $\bu^*(\bx_{f})$ poses a significant challenge for neural network training on the control.  A promising method to overcome this is to exploit a special network structure, the QRnet \cite{QRnet}.
QRnet exploits the solution corresponding to the linear quadratic regulator problem at equilibrium and thus improves the network performance near the equilibrium by providing a good control surrogate.
Suppose we have $\bm u^{\LQR}$, the LQR control for the problem with linearized dynamics and quadratized costs at equilibrium $(\bm x_f, \bm u_f)$. The QRnet can be formulated as
\begin{equation}\label{Eq_QR}
    \bm u^{\QR}(\bm x) = \sigma (\bm u^{\LQR}( \bm x) + \bm u^{\NN}(\bm x;\theta) - \bm u^{\NN}(\bm x_f;\theta)), 
\end{equation}
where $\bm u^{\NN}(\bm x;\theta) $ is any neural network with trainable parameters $\theta$, and  $\sigma$ is a saturating function that satisfies $\sigma(\bm u_f) = \bm u_f$ and $\sigma_{\bm u}(\bm u_f)$ is the identity matrix. Detailed settings are in Appendix~\ref{app:details}.

It is anticipated that $\bu^{\LQR}$ will closely approximate the optimal control near the terminal state $\bx_f$. However, the direct application of QRnet is impractical in our problem. Even with linear and quadratic approximations of $\bm f$ and $\bm L$ around $(\bx_f,\bm u_f)$, the resulting free terminal time problem does not align with the linear quadratic control problem framework, see, \eg~\cite{verriest1991linear} for a detailed discussion. Since minimizing running time is also an objective in the formulation \eqref{eq:cost}, the problem we are addressing is not an equilibrium problem, \ie~$L(\bx_f, \bu^*(t_f)) > 0$.
It should be noted that $\min{L(\bm x_f,\bm u)} = 0$ is valid for stabilization problems \cite{Adaptive_BVP_HJB,QRnet,zang2022machine}, which can be approximated by an LQR problem.
These factors collectively highlight the challenges in applying QRnet to our problem.

To tackle the discontinuity and boost performance near the terminal state, we note a crucial observation. As inferred from \eqref{eq:free_open_problem2}, the open-loop solution for the free terminal time problem is exactly the same as the open-loop solution for the fixed terminal time problem, provided the terminal time is equal to the optimal terminal time. Based on this understanding, we proceed to formulate the fixed terminal time linear quadratic problem in the following manner:
\begin{equation}\label{eq:LQR_problem}
    \begin{aligned}
    & \underset{\bm{x}(t),\bm{u}(t)}{\operatorname{minimize}} \int_{0}^T\tilde{L}(\bm x, \bm u) \rmd t + r_f\| \bm x(t_f) - \bm x_f\|^2 \\
    & s.t.\,  \dot{\bm x}(t) = \tilde{\bm f}(\bm x(t), \bm u(t)), \; \bm x(0) = \bx_0,
\end{aligned}
\end{equation}
where $\tilde{\bm f}$ and $\tilde{\bm L}$ represent the linear and quadratic approximations of $\bm f$ and $\bm L$ at $(\bx_f,\bu_f)$, respectively.
By the properties of LQR, the LQR $\bm \tilde{\bm u}(t,\bm x;T) = \bm u_f+ \bm k(t;T) + \bm K(t;T)(\bm x - \bm x_f)$, where $(\bm k,\bm K)(t;T):[0,T]\rightarrow \bR^7\times \bR^{7\times14}$ satisfies the Riccati equation with terminal time $T$ \cite{sontag2013mathematical}. Due to the time-homogeneity of the LQR problem \eqref{eq:LQR_problem}, $(\bm k,\bm K)(t;T)$ depends only on the remaining time $T-t$. This implies the existence of $(\bm k,\bm K)(t):[0,+\infty) \rightarrow \bR^7\times\bR^{7\times14}$ such that $(\bm k,\bm K)(t;T) = (\bm k,\bm K)(T-t)$. We replace the LQR solution in the original QRnet \eqref{Eq_QR} by
\begin{equation}
 \label{eq:u_lqr}
 \begin{aligned}
         \bm{u}^{\mathrm{LQR}}(\bm x) &=  \bm{\tilde{u}}(0,\bm x;t_f^*(\bm x)) \\
         &= \bm u_f + \bm k(t_f^* (\bm x)) + \bm K(t_f^* (\bm x))(\bm x - \bm x_f),
 \end{aligned}
\end{equation}
which effectively captures the variations of $\bu^*$ to a great extent, as observed on the right of Figure~\ref{fig:terminal_control_lqr}.
Note that in the training phase, we have $t_f^*(\bx^*)$ in \eqref{eq:u_lqr} from the dataset $\mathcal{D}$ \eqref{eq:dataset}, in the actual deployment of closed-loop control, we do not know $t_f^*(\bx^*)$ a priori. To solve this, we train another neural network \(\hat{t}_f(\bm x)\) for predicting optimal terminal time through supervised learning, minimizing the loss function:
\begin{equation}
\label{eq:t_loss}
\underset{\omega}{\operatorname{minimize}}\quad \frac{1}{|\mathcal{D}|} \sum_{i=1}^{|\mathcal{D}|}\left|\hat{t}_f(\bm{x}^{(i)} ; {\omega}) - {t}_f^{(i)} \right|^2.
\end{equation}

\subsection{IVP Enhanced Sampling with Automatic Resampling Times}\label{sec:IVP}

Training closed-loop controllers with a static dataset often fails to yield satisfactory results, an issue attributed to the distribution mismatch phenomenon~\cite{IVP_Enhanced, long2022perturbational}. This phenomenon describes the accumulated error between the state $\bm{x}$ driven by the optimal control $\bm{u}^*$ and a neural network control $\hat{\bm{u}}$, leading to an increasing discrepancy between the distribution of the input states of  $\bm{u}^*$ and $\hat{\bm{u}}$ over time. Consequently, the training data does not accurately represent the states encountered during actual deployment, and the error between $\bm{u}^*$ and $\hat{\bm{u}}$ tends to be larger as time progresses. Our prior research has successfully applied the IVP enhanced sampling~\cite{IVP_Enhanced} approach to address this problem, demonstrating robust performance in both quadrotor and manipulator applications. Nevertheless, the algorithm was primarily designed for tasks with fixed terminal time and necessitated a sequence of predefined times for enhanced sampling, which was heavily reliant on heuristic judgment. In response to these limitations, we have improved the IVP enhanced sampling framework to incorporate automatic resampling times, dubbed IVP-ART, to better encourage exploration. IVP-ART actively modulates the time grid based on the deviation of the state from the optimal path.

\begin{table*}
\centering
\caption{Comparison of success rates and cost ratios among various sampling strategies and iterations. Numbers 0 to 6 correspond to controllers $\hat{\bu}^0$ to $\hat{\bu}^6$, see Section \ref{sec:IVP}. The best performances \revise{after iteration 1} are highlighted in bold.}
\begin{tabular}{cc||ccccccc||c}
\toprule
\textbf{Metrics}                 & \textbf{Sampling}     & \textbf{0} & \textbf{1}           & \textbf{2}           & \textbf{3}           & \textbf{4}           & \textbf{5}           & \textbf{6}           & \textbf{Ensemble}    \\ \midrule
\multirow{2}{*}{\textbf{Success Rate} $\uparrow$} & {\textbf{IVP-ART}}       & \tabincell{c}{0.20 \\ \scriptsize{$\pm$0.12}}  & \tabincell{c}{\textbf{0.63} \\ \scriptsize{$\pm$0.15}}       & \tabincell{c}{\textbf{0.80} \\ \scriptsize{$\pm$0.09}}                 & \tabincell{c}{\textbf{0.90} \\ \scriptsize{$\pm$0.03}}                 & \tabincell{c}{\textbf{0.90} \\ \scriptsize{$\pm$0.05}}                 & \tabincell{c}{\textbf{0.91} \\ \scriptsize{$\pm$0.03}}                 & \tabincell{c}{\textbf{0.92} \\ \scriptsize{$\pm$0.04}}                 & \tabincell{c}{\textbf{0.98} \\ \scriptsize{$\pm$0.00}}                                \\
                                  & \textbf{\DAgger} &  \tabincell{c}{0.27 \\ \scriptsize{$\pm$0.07}}  & \tabincell{c}{0.44 \\ \scriptsize{$\pm$0.21}}       & \tabincell{c}{0.45 \\ \scriptsize{$\pm$0.08}}                 & \tabincell{c}{0.61 \\ \scriptsize{$\pm$0.17}}                 & \tabincell{c}{0.56 \\ \scriptsize{$\pm$0.10}}                 & \tabincell{c}{0.61 \\ \scriptsize{$\pm$0.17}}                 & \tabincell{c}{0.63 \\ \scriptsize{$\pm$0.15}}                 & \tabincell{c}{0.81 \\ \scriptsize{$\pm$0.14}}                                \\  \midrule
\multirow{2}{*}{\textbf{Cost Ratio} $\downarrow$}  & {\textbf{IVP-ART}}     & \tabincell{c}{8.42 \\ \scriptsize{$\pm$0.96}}  & \tabincell{c}{\textbf{4.38} \\ \scriptsize{$\pm$1.37}}       & \tabincell{c}{\textbf{2.91} \\ \scriptsize{$\pm$0.76}}                 & \tabincell{c}{\textbf{2.02} \\ \scriptsize{$\pm$0.20}}                 & \tabincell{c}{\textbf{1.98} \\ \scriptsize{$\pm$0.42}}                 & \tabincell{c}{\textbf{1.95} \\ \scriptsize{$\pm$0.29}}                 & \tabincell{c}{\textbf{1.80} \\ \scriptsize{$\pm$0.36}}                 & \tabincell{c}{\textbf{1.17} \\ \scriptsize{$\pm$0.02}}                                \\
                                  & \textbf{\DAgger}   & \tabincell{c}{7.70 \\ \scriptsize{$\pm$0.61}}       & \tabincell{c}{6.56 \\ \scriptsize{$\pm$1.70}}                 & \tabincell{c}{6.26 \\ \scriptsize{$\pm$0.69}}                 & \tabincell{c}{5.06 \\ \scriptsize{$\pm$1.48}}                 & \tabincell{c}{5.42 \\ \scriptsize{$\pm$0.92}}                 & \tabincell{c}{5.10 \\ \scriptsize{$\pm$1.49}}                 & \tabincell{c}{4.98 \\ \scriptsize{$\pm$1.51}}                 & \tabincell{c}{3.19 \\ \scriptsize{$\pm$1.30}}                 \\ \bottomrule
\end{tabular}
\label{tab: main_sr}
\end{table*}

\begin{figure*}[!ht]
    \centering
    \includegraphics[width=1.0\textwidth]{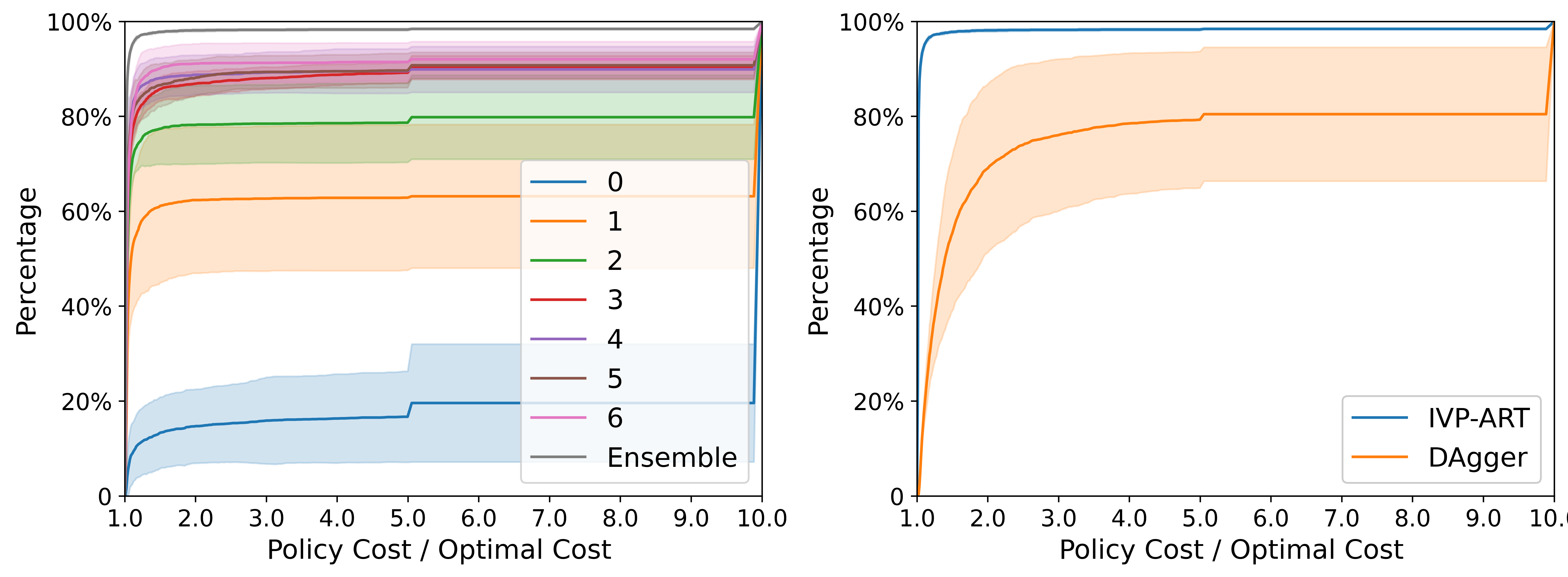}
    \caption{Cumulative distribution of cost ratios (with the ideal curve being a straight horizontal segment passing ratio = 1, percentage=100\%).
    The left panel displays the performances of each iteration of IVP-ART. The right panel compares the performance of the ensemble controller for two sampling strategies.}
    \label{fig:cmp}
\end{figure*}

Let $\mu^0$ be the distribution of initial states of interest.
The IVP-ART algorithm initiates with the training dataset $\mathcal{D}^0$, comprising free terminal open-loop solutions starting from initial states in $X^0$, a set of $N$ initial points independently sampled from $\mu^0$.
A control policy $\hat{\bu}^0$ is then trained on $\mathcal{D}^0$ according to supervision loss \eqref{eq:u_loss}.
At $k$-th iteration of IVP-ART, we compute IVP trajectories 
for states in $X^0$ using the latest closed-loop control policy $\hat{\bu}^{k-1}$, as specified by \eqref{eq:IVP}.
For each trajectory, the first state that deviates from its optimal counterpart by a margin greater than $\tau>0$ is identified and recorded:
\begin{equation}\label{eq:t_ivp}
\begin{aligned}
\bx^{\hat{\bu}}(\bx_0;\tau)& = \bm{x}^{\bm{\hat{u}}}(t^{\bm{\hat{u}}}(\bm x_0; \tau); \bx_0), \\
    t^{\bm{\hat{u}}}(\bm x_0; \tau)
    & = \inf \{t:\| \bm{x}^{\bm{\hat{u}}}(t; \bm x_0) - \bm{x}^*(t; \bm x_0) \| > \tau \}.
\end{aligned}
\end{equation}
Here $t^{\hat{\bu}}(\bx_0;\tau)$ is the automatic resampling time, adaptively generated for each IVP trajectory.
Subsequently, we define $X^{k}=\{\bx^{\hat{\bu}^{k-1}}(\bx_0;\tau):\bx_0\in X^0\}$ and construct optimal open-loop solutions starting from states in $X^k$ to form the dataset $\hat{\mathcal{D}}^k$.
A new controller $\hat{\bu}^k$ is trained starting with randomly initialized weights on the combined dataset $\mathcal{D}^k = \mathcal{D}^{k-1} \cup \hat{\mathcal{D}}^k$.
We shall remark that when employing the modified QRnet architecture introduced in Section \ref{sec:qrnet}, each control network is paired with a terminal time network adhering to the supervision loss described in \eqref{eq:t_loss}.
Both networks are trained on the same dataset $\mathcal{D}^k$.

This iterative process is conducted through $K$ cycles for $k=1, \ldots, K$. Upon completion, we obtain a series of $K+1$ networks: $\hat{\bm{u}}^0, \hat{\bm{u}}^1, \ldots, \hat{\bm{u}}^{K}$. In the final step, we construct an ensemble controller by averaging the individual controllers (excluding $\hat{\bu}^0$):
\(
\hat{\bm{u}}^{\text{Ensemble}} = \sum_{k=1}^{K} \hat{\bm{u}}^k/K
\).
According to the numerical experiments detailed in Section \ref{sec:experiment_results}, this ensemble controller demonstrates superior performance, outperforming controllers in all previous iterations significantly.
For a comprehensive overview of the IVP-ART algorithm, readers are referred to Algorithm~2 in Appendix~\ref{Sec:Alg}.

\section{Experiments}

\subsection{Experimental Settings}
Comprehensive details regarding our experimental setup are detailed in Appendix~\ref{app:details}.

\paragraph{Data Generation} It is paramount to highlight that the domain for initial states in our study is exceptionally broad, substantially surpassing the scope utilized in prior research~\cite{IVP_Enhanced}. Specifically, the initial positions $\bm{q}_0$ are uniformly and independently sampled from an expansive 7-dimensional cube centered around $\bm{q}_c$ with side length $1.0$ while the initial velocities $\bm{v}_0$ are set to zero. 
The initial domain is $50^7 \approx 7.8 \times 10^{11}$ times larger than the cubic domain with a side length of $0.02$, as used in \cite{IVP_Enhanced}, significantly extending the scope towards a more realistic problem setting.
Leveraging the iterative DDP-based free terminal time open-loop optimal control solver described in Section \ref{sec:solver}, we manage to generate datasets comprising 500 training trajectories for initial training, 100 validation trajectories, and 1000 testing trajectories, achieving the convergence rate of 96.7\%, a notable improvement over the rate of 82.8\% without proposed marching of time discretization.

\paragraph{Evaluation} 
In evaluation, we conduct IVP simulations~\eqref{eq:IVP}, extending the time span to $2.0$ seconds. This duration surpasses the maximal optimal terminal time ($\approx 1.4 $ seconds) in our dataset (see Figure~\ref{fig:optimal_time} in Appendix~\ref{app:details}).
To assess closed-loop control performance, we employ two metrics.
The first metric is the success rate, which measures the ability of the control to guide the manipulator to the terminal states.
A terminal state is deemed reached if the manipulator is within an $\epsilon = 0.001$ Euclidean distance, \ie\, $\|\bm{x} - \bm{x}_f\| \leq \epsilon$.
The second metric evaluates the optimality of the closed-loop control through cost ratios. 
This ratio compares the total cost of simulated trajectories to the optimal cost. 
To address extreme values, we cap the ratio at 10.0 for trajectories that fail to reach the terminal state and at 5.0 for those that succeed but incur excessive costs.

\subsection{Evaluation Results}\label{sec:experiment_results}

In our evaluation of network architectures, we found that the Multi-Layer Perceptron (MLP) controller performed very poorly, with a success rate of zero due to its inability to handle discontinuities at the terminal state, as discussed in Section~\ref{sec:qrnet}. See also Figure~\ref{fig:lqr_traj} in Appendix~\ref{app:results}. This highlights the effectiveness of our enhanced QRnet in addressing discontinuities in free terminal time problems. Subsequent experiments exclusively utilized the QRnet architecture. \revise{The LQR solution at the terminal state is computed by PyDrake~\cite{drake}.}
We also compared our approach with the \DAgger\ (Dataset Aggregation) approach~\cite{dagger}, another adaptive sampling method to address distribution mismatches, with detailed experimental setups in Appendix~\ref{app:details}. The mean and standard deviation of the evaluation metrics, based on three independent trials, are reported in Table~\ref{tab: main_sr}.
The initial performance of the control network is comparable across methods since the initial training dataset is the same.
Despite similar initial performances, IVP-ART significantly outperformed \DAgger\ in cost ratios and success rates from the first iteration, with both methods reaching a performance plateau.
Despite similar initial performances, IVP-ART significantly outperforms \DAgger\ in cost ratios and success rates from the first iteration, even though both methods exhibit consistent improvements throughout the iterative process.
It is worth noting that both methods achieve a steady state performance, reaching a plateau where further improvements are insignificant.
The ensemble network for IVP-ART demonstrates exceptional performance, achieving a success rate of 98\% and a cost ratio of 1.17, which is remarkably close to the theoretical optimum of 1.0.
Additionally, we visualize the cumulative distribution of cost ratios in Figure~\ref{fig:cmp}.
Overall, the results demonstrate that the integration of IVP-ART with a modified QRnet architecture not only outperforms other strategies but also closely aligns with the optimal policy.
Exceptionally, the IVP-ART ensemble network achieved a 98\% success rate and a cost ratio of 1.17, close to the theoretical optimum of 1.0. These results confirm that integrating IVP-ART with the modified QRnet architecture not only surpasses other strategies but also closely approximates the optimal policy.

\section{Conclusion}

In this study, we have advanced the supervised learning approach to create a free terminal time optimal controller for robotic manipulators. This involved generating a dataset from optimal open-loop trajectories using numerical solvers and training a policy network using that data. Key advancements include a marching scheme in the open-loop solver for better solutions, adapting the QRnet for enhanced stability at terminal states, and an innovative adaptive IVP sampling method for improved dataset quality. These combined efforts led to a versatile closed-loop policy capable of handling diverse tasks with near globally optimal total costs. In the future, we plan to incorporate constraints in the manipulator problem, including the avoidance of self-collision and collisions with the environment. Additionally, we aim to develop a closed-loop controller capable of handling varying terminal states.

\bibliographystyle{IEEEtran} 
\bibliography{IEEEabrv,ref}

\newpage
\appendices
\section{Algorithm}\label{Sec:Alg}
\begin{algorithm}[!ht]
\caption{Iterative DDP-based Free Terminal Time Open-loop Optimal Control Solver}
\label{alg:solver}
\SetAlgoLined
\SetKwFunction{FMain}{Main}
\SetKwFunction{FOpenSolver}{OpenLoopSolver}

\DontPrintSemicolon
\SetKwProg{Fn}{Function}{:}{\KwRet{$\bu^*(t;\bx_0, t_f^*)$, $\bx^*(t;\bx_0, t_f^*)$, $t_f^*$}}
\Fn{\FMain{$\cdots$}}{
    \KwData{Initial state $\bx_0$, initial terminal time estimate $t_f^0$, sequence of time step counts $\{n_1, n_2, \cdots, n_J\}$, 
    terminal time maximal update fraction $\alpha$,
    convergence threshold $\epsilon$,
    threshold to switch from gradient descent to quasi-Newton method $\epsilon_0$,
    unified time step size $\Delta t$}
    $k\gets 0$\;
    $gd\gets True$\;
    
    Set $\bm u^{-1}_J=0$ and $\bm x^{-1}_J=0$

    \While{True}{
        $k \gets k+1$\;
    
        Set $\bm u^k_0=u^{k-1}_J$ and $\bm x^{k}_0=\bm x^{k-1}_J$
    
        $\bm u_{J}^k, \bm x_{J}^k \gets $ \FOpenSolver{$\bx_0$, $t_f^k$, $\{n_1, n_2, \cdots, n_J\}$, $\{\bu_0^k, \bx_0^k\}$}
    
        Compute gradient $\rmd t_f^{k}$ at $t_f^{k}$ according to \eqref{eq:graident_t} with $\bm u_{J}^k$ and $\bm x_{J}^k$\;
                
        \lIf{$|\rmd t_f^{k}|<\epsilon$}{break}
        \eIf{$gd$}{
            Set $\delta t_f^{k} = \rmd t_f^{k}$\;
            
            \lIf{$|\rmd t_f^{k}|<\epsilon_0$}{$gd\gets False$}\;
        }{
            Compute $\delta t_f^{k}$ using quasi-Newton method: $\delta t_f^{k} = \frac{t_f^{k} - t_f^{k-1}}{\rmd t_f^{k} - \rmd t_f^{k-1}}\rmd t_f^{k}$.
        }
        $\alpha^k \gets \alpha|t_f^{k}|/| \delta t_f^{k}|$\;

        Update $t_f^{k+1} = t_f^{k} - \alpha^{k} \delta t_f^{k}$\;
    }
    $N \gets \text{Round}(t_f^k/\Delta t)$\;

    $t_f^*\gets N\Delta t$\;

    $\bu^*, \bx^* \gets $ \FOpenSolver{$\bx_0$, $t_f^*$, $\{n_1, n_2, \cdots, n_J, N\}$, $\{\bu_J^k, \bx_J^k\}$}
}
\SetKwProg{Fn}{Function}{:}{\KwRet{$\bu_J,\bx_J$}}
\Fn{\FOpenSolver{$\cdots$}}{
    \tcp{Robust Open-Loop Solver}
    \KwData{Initial state $\bx_0$, terminal time $t_f$, sequence of time step counts ${n_1, n_2, \cdots, n_J}$, initial guess $\bu_0, \bx_0$}
    \For{$j \leftarrow 1$ \KwTo $J$}
    {
        Solve fixed terminal time open-loop problem with terminal time $t_f$ and $n_j$ time steps using $\bm u_{j-1}$ and $\bm x_{j-1}$ as initial guess to obtain $\bm u_{j}$ and $\bm x_{j}$\;
    }
}
\end{algorithm}

\begin{algorithm}[!htpb]
\caption{IVP Enhanced Sampling with Automatic Resampling Times}
\label{alg:sampling}
\SetAlgoLined
\DontPrintSemicolon
\KwIn{Initial distribution $\mu_0$, number of initial points $N$, number of iterations $K$, state error margin $\tau$.}

\textbf{Solver:} The open-loop optimal control solver $\mathcal{S}$ for free terminal time problems, \eg  \ Solver~\ref{alg:solver}.

\tcp{Initialization and Training}
Independently sample $N$ initial points from $\mu_0$ to get an initial point set $X_0$.

Call the solver $\mathcal{S}$ to obtain open-loop solutions as $\mathcal{D}_0 = \{(t,\bm x^*(t;\bm x_0),\bm u^*(t;\bm x_0)): \bm x_0 \in {X}_0\}$.

Train $\hat{\bm u}_{0}$ with the initial dataset $\mathcal{D}_0$.

\tcp{Training with Adaptive Samples}
\For{$i=1$ \KwTo $K$}{
    For any $\bm x_0 \in  X_0$, compute IVP $(\bm x_0, \hat{\bm u}^{i-1})$ according to \eqref{eq:dyanmic_system}.  \Comment{Exploration}

    Generate a set of new states $X^{k}=\{\bx^{\hat{\bu}^{k-1}}(\bx_0;\tau):\bx_0\in X^0\}$ according to~\eqref{eq:t_ivp} with margin $\tau$.

    For any $\bm x_i \in {X}_i$, call the solver $\mathcal{S}$ to obtain $\bm x^*(t;\bm x_i) $ and $\bm u^*(t;\bm x_i)$ and to form   $\hat{\mathcal{D}}^i$ \Comment{Labeling}
    
    Set $\mathcal{D}^i= \mathcal{D}^{i-1} \cup \hat{\mathcal{D}}^i$.
    
    Train $\hat{\bm u}_{i}$ on $\mathcal{D}^i$. \Comment{Training}
}

\KwOut{$\{\hat{\bm u}^{0}, \dots, \hat{\bm u}^{K}\}$ and ensemble network $\sum_{k=1}^K\hat{\bm u}^{k}/K$.}

\end{algorithm}

\section{Additional Experiment Details}\label{app:details}

\paragraph{Manipulator.} We use similar configuration as~\cite{IVP_Enhanced} with some adaption for the free terminal time setting.
We set the center of initial states as $ \bm{q}_c=[1.6800,1.2501,2.4428,$ $-1.2669,-0.9778,1.1236,-1.3575]^T$ and the terminal state as $\bm{q}_f=[2.7736,0.5842,1.5413,-1.7028,-2.1665,0.0847,$ $-2.5764]^T$.
The parameters in the running cost~\eqref{eq:running_cost_used} are $r_t=100, r_a=0.005, r_u=0.025$, and the $r_f$ for the terminal cost in Section \ref{sec:solver} is $2.5e5$.
To compute the total cost for the trajectories, we adhere to a criterion where the terminal time is defined as the moment when the distance of the state to the terminal state falls within $\epsilon = 0.001$. At this point, we integrate the running cost up to the terminal time without terminal cost, resulting in the total cost for the trajectory.

\paragraph{Free Time Solver Hyperparameters.}
We present the hyperparameters used in Algorithm \ref{alg:solver}.
The sequence of time step counts is $\{150, 300, 450, 600, 750, 900, $ $ 1200, 1500, 1750\}$.
Other parameters are summarized in Table \ref{tab:hyperparameters}.

 \begin{table}[h]
\centering
\begin{tabular}{|c|c|c|}
\hline
Hyperparameter & Value & Description \\
\hline
$t_f^0$ & 1.2 & initial estimation of terminal time \\
$\alpha$ & $0.2$ & terminal time maximal update fraction \\
$\epsilon$ & $1e-6$ & convergence threshold \\
$\epsilon_0$ & $0.3$ & threshold to switch to quasi-Newton method \\
$\Delta t$ & $5e-4$ & unified time step size \\
\hline
\end{tabular}
\caption{Summary of hyperparameters in the free time solver Algorithm \ref{alg:solver}}
\label{tab:hyperparameters}
\end{table}

\paragraph{Additional Data Generation during Resampling} When generating $\hat{\mathcal{D}}^k$ with initial states from $X^k$ to perform the IVP-ART algorithm in Section \ref{sec:IVP}, we set the initial guess of the terminal time for $\bx^{\hat{\bu}^{k-1}}(\bx_0; \tau)$ as $t_f^*(\bx_0) - t^{\hat{\bu}^{k-1}}(\bx_0;\tau)$ and use the original optimal trajectories, $(\bx^*(t; \bx_0),\bu^*(t; \bx_0))$ for $t \ge t^{\hat{\bu}^{k-1}}(\bx_0;\tau)$ as the initial guess for the fixed time problem, where $t^{\hat{\bu}^{k-1}}(\bx_0;\tau))$ is the automatic resampling time in \eqref{eq:t_ivp}.
The convergence rate when generating $\hat{\mathcal{D}}^k$ is very high, above 99\% in all iterations.

\paragraph{Network Architectures.}
We utilized similar architectures for the control network and the terminal time network. Specifically, the control network is implemented in two variants: QRnet and the vanilla MLP. The internal architecture of the QRnet, $\bu^{\NN}$ in \eqref{Eq_QR}, is identical to that of the vanilla MLP. The control MLP is designed with an input dimension of 14 and comprises 4 hidden layers, with layer sizes of 32, 64, 64, and 32, respectively. The output dimension for the control network is 7. For activation functions, the first two layers utilize the tanh function, while the remaining layers employ the Exponential Linear Unit (ELU).

Parallel to the control network, the terminal time network features an identical structure in terms of the number of hidden layers and neuron sizes. However, it exclusively uses ELU as the activation function in all hidden layers. The output dimension of this network is 1. To ensure the positiveness of the output, a softplus activation function is added at the output layer.

\paragraph{Design for QRnet.}

The saturation function $\sigma(\cdot)$ used in~\eqref{Eq_QR} is defined on a coordinate-wisely basis:
\begin{equation*}
    \sigma (u) = u_\tmin + \frac{u_\tmax - u_\tmin}{1+c_1 \exp[-c_2(u-u_1)]},
\end{equation*}
where 
$ c_1=(u_\tmax - u_1) / (u_1-u_\tmin), c_2=(u_\tmax - u_\tmin) / [(u_\tmax -u_1)(u_1-u_\tmin)]$ with
$u_\tmin, u_\tmax$ being minimum and maximum values for $u$. Here $u, u_{\min}=-2000$ and $u_{\max}=2000$ are the corresponding values at each coordinate of $\bm u$, $ \bm u_\text{min}$, $ \bm u_\text{max}$, respectively.
Regarding $\bm{u}^\text{LQR}$ in QRnet~\eqref{Eq_QR}, its effectiveness is predominantly near the terminal state $\bm{x}_f$. Consequently, the linear quadratic regulator is not required for the entire duration.
Therefore, rather than applying $\bm{k}, \bm{K}$ from Equation~\eqref{eq:u_lqr} for all $t \ge 0$, these are employed only for a limited duration ($t < t_m$) and replaced with $\bm{u}_f$ for larger time frames ($t > t_M$).
Specifically, we adopt
\[
\tilde {\bm K} = s\bm K + (1-s) 0, \tilde {\bm k} = s\bm k + (1-s) 0
\]
where $s = s(t)$ is defined as
\[
\resizebox{0.95\linewidth}{!}{$
s = s(t) = 
\begin{cases} 
1 & \text{if } t < t_m, \\
\text{Sigmoid}\left(-(t - t_m) \dfrac{t_\tmax - t_\tmin}{t_M - t_m} + t_\tmax\right) & t_m \le t \le t_M,\\
0 & \text{if } t > t_M.
\end{cases}
$}
\]
In this setup, $t_{\text{min}}$ and $t_{\text{max}}$ are selected such that $\text{Sigmoid}(t_{\text{min}}) \le \varepsilon$ and $\text{Sigmoid}(t_{\text{max}}) \ge 1 - \varepsilon$.
In the experiments, we solve the LQR \eqref{eq:LQR_problem} with a duration $T=0.8$ seconds, setting $t_m=0.08$ seconds, $t_M=T=0.8$ seconds, and $\varepsilon=1e-5$. These parameters are flexible and do not alter the essence of the problem.

\paragraph{Training Details.} In our study, we employ an iterative training approach, where each iteration consists of 200 epochs. We carry out a total of 6 iterations ($K=6$) for network training. The batch size for our training process is set at 1024. We use Adam~\cite{kingma2015adam} with a learning rate of 0.001. To ensure the robustness of our model, we incorporate a validation step every 10 epochs. The model that demonstrates the best performance on the validation set is saved as the representative model for that particular iteration. This iterative training and validation process ensures that our model is both accurate and reliable. 

 We set the state error threshold as $\tau=1.0$ in IVP-ART for experiments in Section~\ref{sec:experiment_results} and conduct additional experiments with $\tau=0.5$ and $2.0$ in Appendix~\ref{app:results}.

The \DAgger~algorithm implementation necessitates a predefined temporal grid for data updates, denoted as $0< t_1 < \dots < t_f$. To accommodate the final time $t_f$ for each trajectory in the free terminal time setting, we establish an individual time grid for each trajectory, setting $t_1=0.25 t_f$ and $t_2=0.75 t_f$. 
The training process of \DAgger\ parallels that of IVP-ART, with a notable distinction: the collection of new states $X^k$ occurs on predetermined grids rather than the adaptive time grids utilized in IVP-ART, as delineated in equation~\eqref{eq:t_ivp}.

\begin{figure}[!ht]
    \centering
    \includegraphics[width=0.475\textwidth]{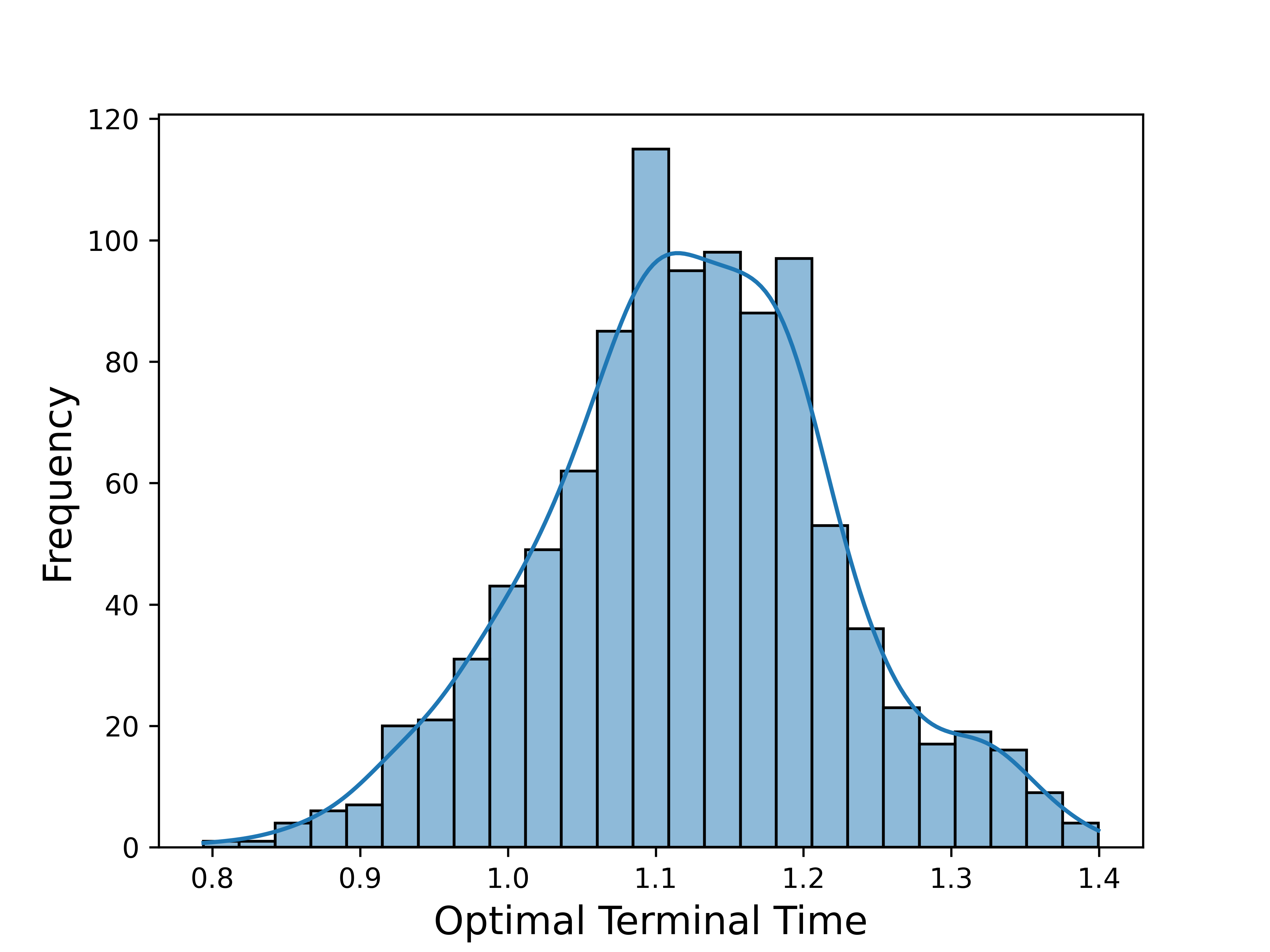}
    \caption{Distribution of optimal terminal time for initial states of the manipulator problem.}
    \label{fig:optimal_time}
\end{figure}

\begin{table*}[!h]
\centering

\begin{tabular}{cc||ccccccc||c}
\toprule
\textbf{Success Rate} $\uparrow$       & \textbf{Threshold} & \textbf{0} & {\textbf{1}} & {\textbf{2}} & \textbf{3}           & {\textbf{4}} & \textbf{5}           & {\textbf{6}} & {\textbf{Ensemble}} \\ \midrule
 & \textbf{$\tau = 0.5$}       & \tabincell{c}{0.30 \\ \scriptsize{$\pm$0.11}}  & \tabincell{c}{{0.53} \\ \scriptsize{$\pm$0.06}}       & \tabincell{c}{{0.70} \\ \scriptsize{$\pm$0.09}}                 & \tabincell{c}{{0.72} \\ \scriptsize{$\pm$0.14}}                 & \tabincell{c}{{0.59} \\ \scriptsize{$\pm$0.18}}                 & \tabincell{c}{{0.75} \\ \scriptsize{$\pm$0.09}}                 & \tabincell{c}{{0.71} \\ \scriptsize{$\pm$0.17}}                 & \tabincell{c}{{0.87} \\ \scriptsize{$\pm$0.06}}                                \\
      \textbf{Union}  & \textbf{$\tau = 1.0$}       & \tabincell{c}{0.20 \\ \scriptsize{$\pm$0.12}}  & \tabincell{c}{{0.63} \\ \scriptsize{$\pm$0.15}}       & \tabincell{c}{{0.80} \\ \scriptsize{$\pm$0.09}}                 & \tabincell{c}{{0.90} \\ \scriptsize{$\pm$0.03}}                 & \tabincell{c}{{0.90} \\ \scriptsize{$\pm$0.05}}                 & \tabincell{c}{{0.91} \\ \scriptsize{$\pm$0.03}}                 & \tabincell{c}{{0.92} \\ \scriptsize{$\pm$0.04}}                 & \tabincell{c}{{0.98} \\ \scriptsize{$\pm$0.00}}                                \\
                                  & \textbf{$\tau = 2.0$} & \tabincell{c}{0.30 \\ \scriptsize{$\pm$0.03}}       & \tabincell{c}{0.90 \\ \scriptsize{$\pm$0.01}}                 & \tabincell{c}{0.93 \\ \scriptsize{$\pm$0.06}}                 & \tabincell{c}{0.94 \\ \scriptsize{$\pm$0.02}}                 & \tabincell{c}{0.96 \\ \scriptsize{$\pm$0.02}}                 & \tabincell{c}{0.91 \\ \scriptsize{$\pm$0.05}}                 & \tabincell{c}{0.98 \\ \scriptsize{$\pm$0.00}}                 & \tabincell{c}{1.00\footnotemark \\ \scriptsize{$\pm$0.00}}                 \\  \midrule
\textbf{Replacement}  & \textbf{$\tau=1.0$}      & \tabincell{c}{0.28 \\ \scriptsize{$\pm$0.09}}  & \tabincell{c}{0.76 \\ \scriptsize{$\pm$0.06}}       & \tabincell{c}{0.86 \\ \scriptsize{$\pm$0.06}}                 & \tabincell{c}{0.76 \\ \scriptsize{$\pm$0.02}}                 & \tabincell{c}{0.91 \\ \scriptsize{$\pm$0.02}}                 & \tabincell{c}{0.91 \\ \scriptsize{$\pm$0.02}}                 & \tabincell{c}{0.78 \\ \scriptsize{$\pm$0.04}}                 & \tabincell{c}{0.98 \\ \scriptsize{$\pm$0.01}}                                \\
 \bottomrule
\end{tabular}
 \caption{Comparison of success rates among various dataset integrations, thresholds and iterations. The setting of ``Union, $\tau=1.0$" is the same as ``IVP-ART" in Table~\ref{tab: main_sr}.}\label{tab: ablation_sr}
\end{table*}

\begin{table*}[!h]
\centering
\begin{tabular}{cc||ccccccc||c}
\toprule
\textbf{Cost Ratio} $\downarrow$    & \textbf{Threshold} & \textbf{0} & {\textbf{1}} & {\textbf{2}} & \textbf{3}           & {\textbf{4}} & \textbf{5}           & {\textbf{6}} & {\textbf{Ensemble}} \\ \midrule
 & \textbf{$\tau = 0.5$}       & \tabincell{c}{7.83 \\ \scriptsize{$\pm$0.54}}  & \tabincell{c}{{5.41} \\ \scriptsize{$\pm$0.52}}       & \tabincell{c}{{3.83} \\ \scriptsize{$\pm$0.83}}                 & \tabincell{c}{{3.64} \\ \scriptsize{$\pm$1.16}}                 & \tabincell{c}{{4.88} \\ \scriptsize{$\pm$1.70}}                 & \tabincell{c}{{3.38} \\ \scriptsize{$\pm$0.83}}                 & \tabincell{c}{{3.70} \\ \scriptsize{$\pm$1.57}}                 & \tabincell{c}{{2.24} \\ \scriptsize{$\pm$0.61}}                                \\
    \textbf{Union}            & \textbf{$\tau = 1.0$} & \tabincell{c}{8.42 \\ \scriptsize{$\pm$0.96}}  & \tabincell{c}{{4.38} \\ \scriptsize{$\pm$1.37}}       & \tabincell{c}{{2.91} \\ \scriptsize{$\pm$0.76}}                 & \tabincell{c}{{2.02} \\ \scriptsize{$\pm$0.20}}                 & \tabincell{c}{{1.98} \\ \scriptsize{$\pm$0.42}}                 & \tabincell{c}{{1.95} \\ \scriptsize{$\pm$0.29}}                 & \tabincell{c}{{1.80} \\ \scriptsize{$\pm$0.36}}                 & \tabincell{c}{{1.17} \\ \scriptsize{$\pm$0.02}}                                \\
                                  & \textbf{$\tau = 2.0$} & \tabincell{c}{7.74 \\ \scriptsize{$\pm$0.45}}       & \tabincell{c}{2.03 \\ \scriptsize{$\pm$0.14}}                 & \tabincell{c}{1.73 \\ \scriptsize{$\pm$0.58}}                 & \tabincell{c}{1.66 \\ \scriptsize{$\pm$0.23}}                 & \tabincell{c}{1.41 \\ \scriptsize{$\pm$0.12}}                 & \tabincell{c}{2.10 \\ \scriptsize{$\pm$0.52}}                 & \tabincell{c}{1.31 \\ \scriptsize{$\pm$0.02}}                 & \tabincell{c}{1.05 \\ \scriptsize{$\pm$0.00}}                 \\  \midrule
\textbf{Replacement}  & \textbf{$\tau=1.0$}      & \tabincell{c}{7.63 \\ \scriptsize{$\pm$0.71}}  & \tabincell{c}{3.21 \\ \scriptsize{$\pm$0.47}}       & \tabincell{c}{2.41 \\ \scriptsize{$\pm$0.57}}                 & \tabincell{c}{3.50 \\ \scriptsize{$\pm$0.32}}                 & \tabincell{c}{1.90 \\ \scriptsize{$\pm$0.20}}                 & \tabincell{c}{2.00 \\ \scriptsize{$\pm$0.30}}                 & \tabincell{c}{3.14 \\ \scriptsize{$\pm$0.38}}                 & \tabincell{c}{1.25 \\ \scriptsize{$\pm$0.10}}                                \\
 \bottomrule
\end{tabular}
\caption{Comparison of mean of cost ratios among various dataset integrations, thresholds and iterations. The setting of ``Union, $\tau=1.0$'' is the same as ``IVP-ART'' in Table~\ref{tab: main_sr}.}\label{tab: ablation_cr}
\end{table*}

\begin{figure}[h]
    \centering
    \includegraphics[width=0.475\textwidth]{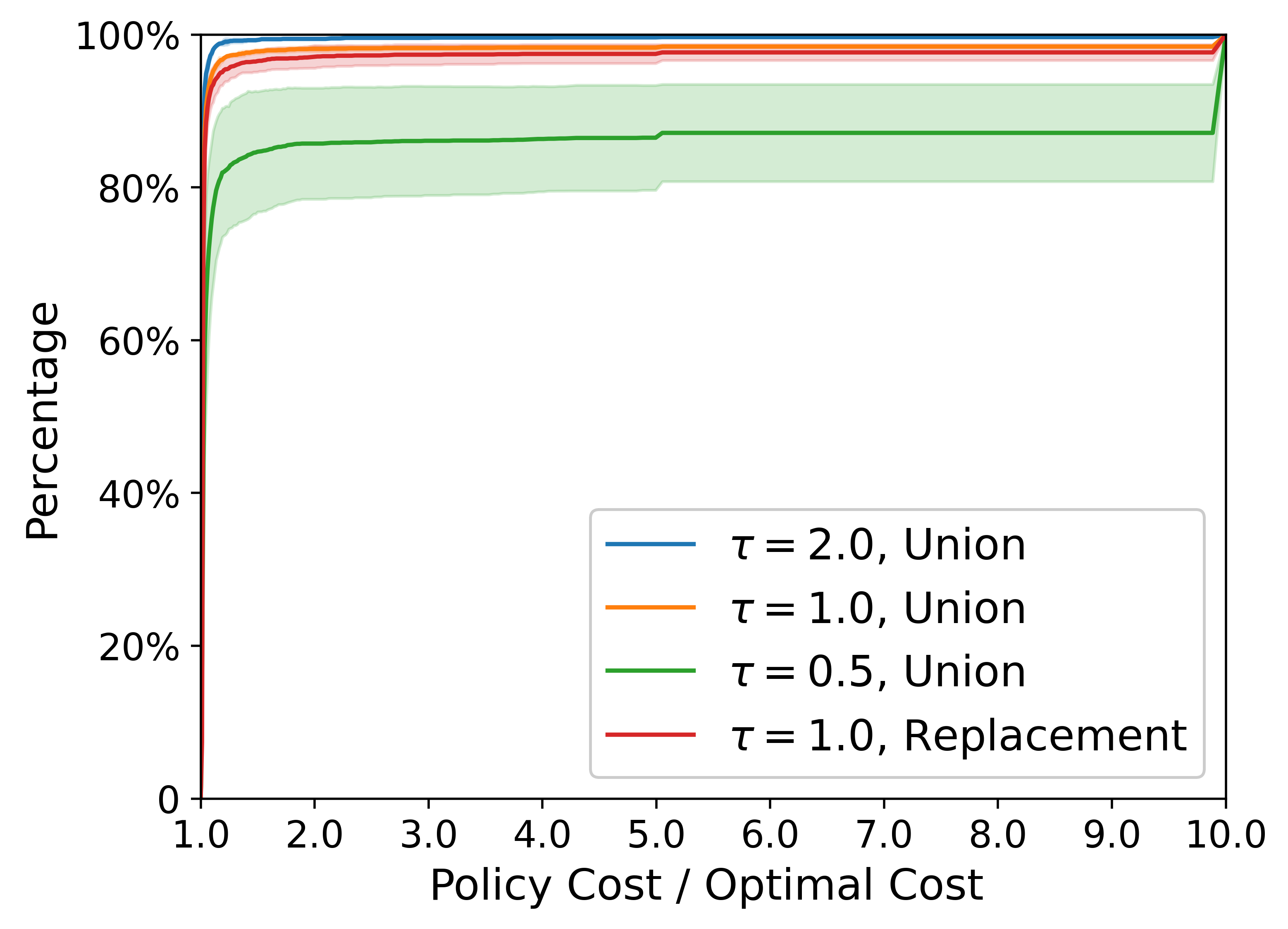}
    \caption{Cumulative distribution of cost ratios (with the ideal curve being a straight horizontal segment passing ratio = 1, percentage=100\%). }
    \label{fig:cmp_app}
\end{figure}

\begin{figure*}[!h]
    \centering
    \includegraphics[width=0.6\textwidth]{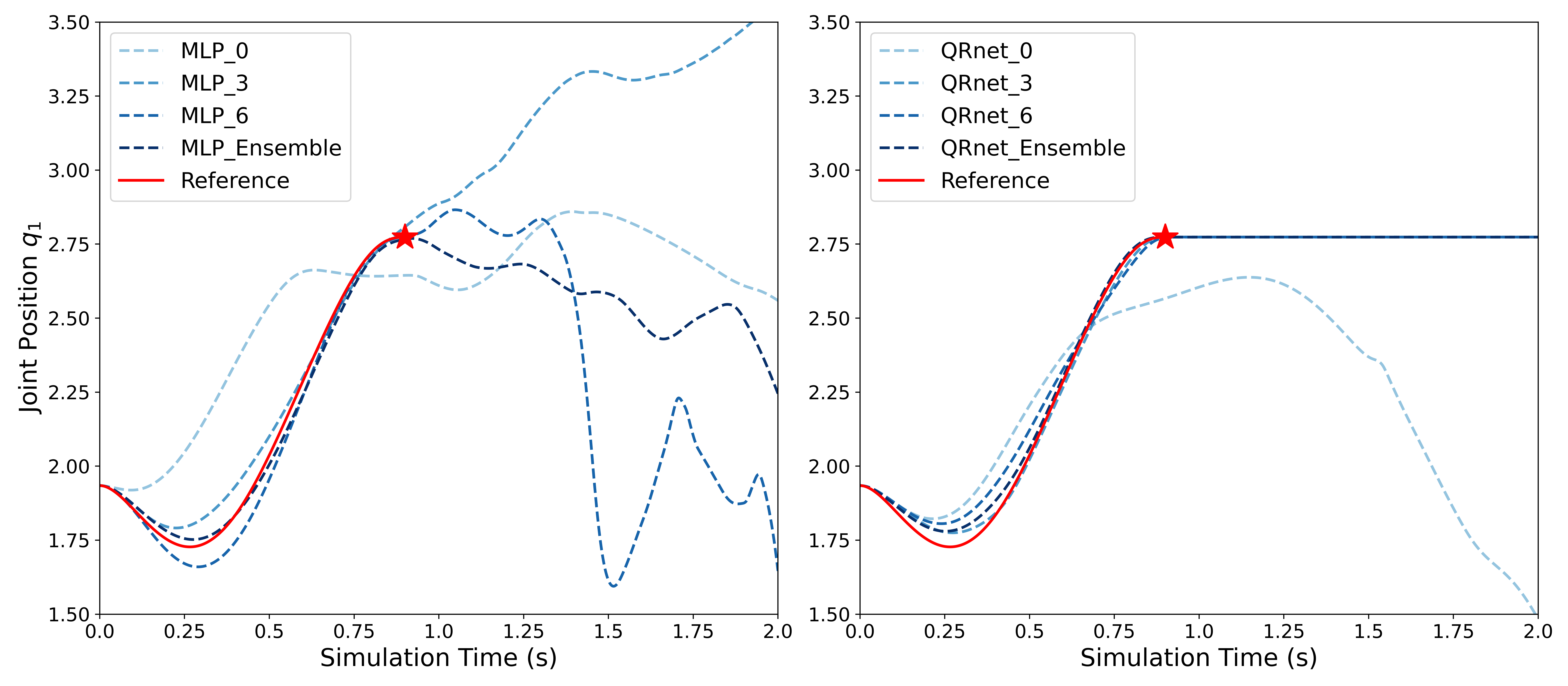}
    \caption{Comparative trajectories of joint position $q_1$. The left panel showcases IVPs governed by MLPs across varying iterations or ensemble, whereas the right panel represents that of the QRnet. The trajectories of the IVPs are drawn in blue, the reference solution, serving as the optimal path, is highlighted in red, and the desired terminal time-state pair is indicated by a red star.}
    \label{fig:lqr_traj}
\end{figure*}

\section{Additional Results}\label{app:results}
In order to substantiate the robustness of our work, we have conducted additional experiments with different dataset integrations as shown in Table~\ref{tab: ablation_sr}, Table~\ref{tab: ablation_cr} and Figure~\ref{fig:cmp_app}. These include an evaluation of state error thresholds at $\tau=0.5$ and $2.0$. The performance observed at these thresholds was found to be comparable to the results presented in Table~\ref{tab: main_sr} for $\tau=1.0$. 
Besides, in the study by \cite{IVP_Enhanced}, new datasets were formed by replacing subsequent trajectories with newly generated ones. Unlike this approach, our IVP-ART method retained these data in the training dataset. However, to enable a more comprehensive comparison, we shifted our IVP-ART approach from utilizing a union of the trajectories to replacing them, in line with the method used in \cite{IVP_Enhanced}. Our adapted method still demonstrates superior performance compared to \DAgger~\cite{dagger}, and it is slightly outperformed by the original version of IVP-ART. This adaptation served as an ablation study, further validating the effectiveness of our IVP-ART method.

To compare MLP and QRnet further, we have plotted the trajectories of the first component of the position $\bm q$ in Figure~\ref{fig:lqr_traj}. The IVP trajectories controlled by MLPs with different iterations ($\hat{\bm{u}}_i, i=0, 3, 6$) or ensemble are presented on the left panel, while the ones controlled by QRnet are displayed on the right panel. Given the same initial state $\bm x_0$ and the reference trajectory, the IVP controlled by an MLP policy might be close to the optimal trajectory, but ultimately, it fails to achieve stability at the terminal state. In contrast, utilizing QRnet successfully maintains stability around the terminal state. 

\footnotetext{The success rate is rounded to two decimal places for presentation purposes. The actual success rate is 0.997.}

\begin{figure*}[h]
    \centering
    \includegraphics[width=0.15\textwidth]{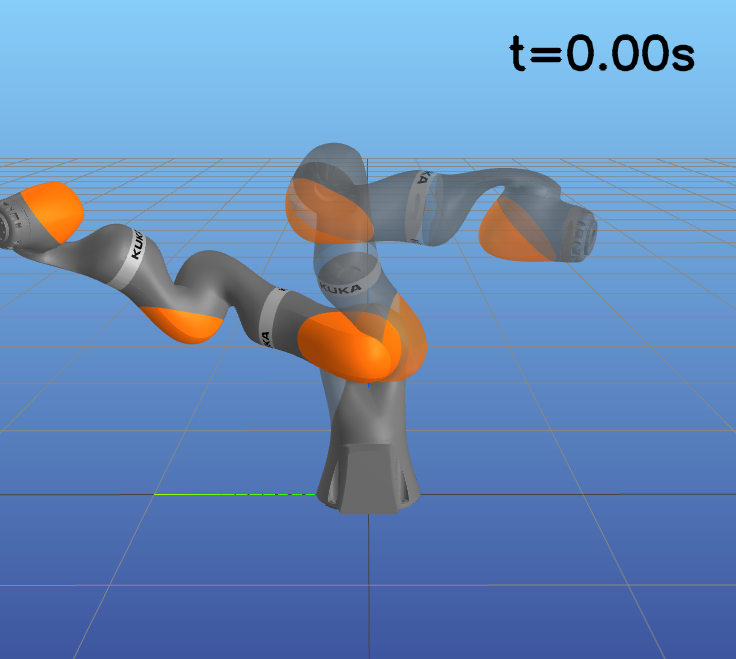}
    \includegraphics[width=0.15\textwidth]{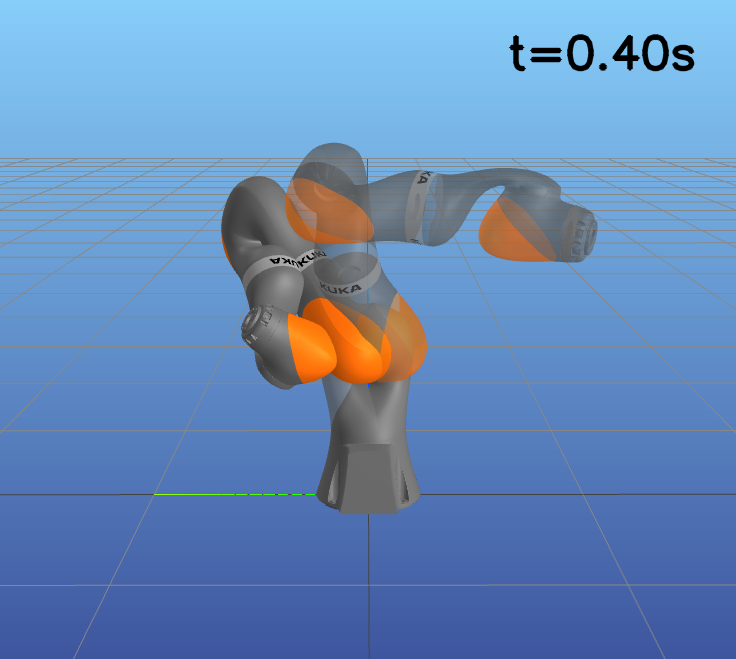}
    \includegraphics[width=0.15\textwidth]{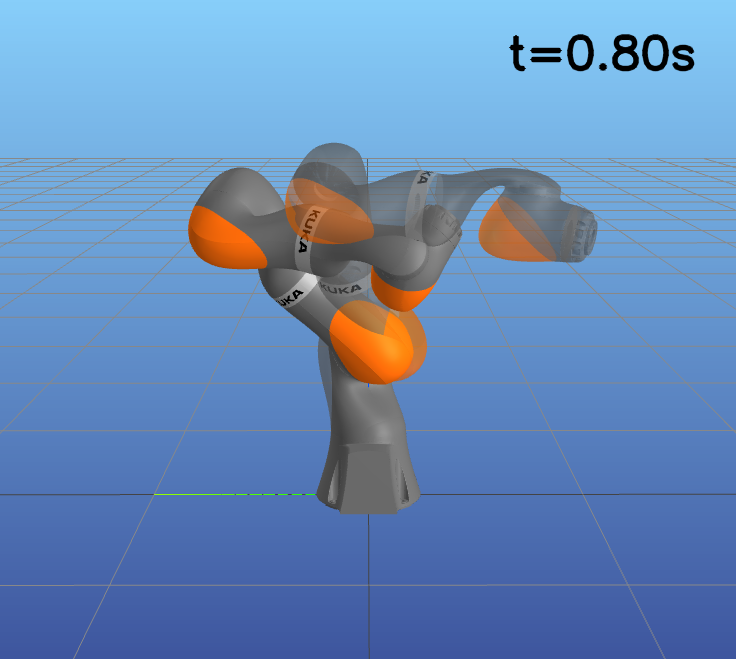}
    \includegraphics[width=0.15\textwidth]{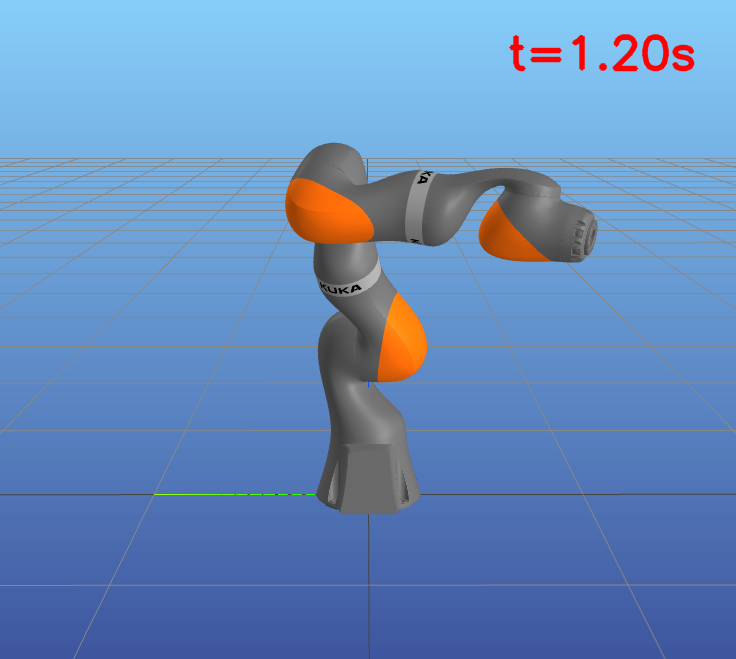}
    \includegraphics[width=0.15\textwidth]{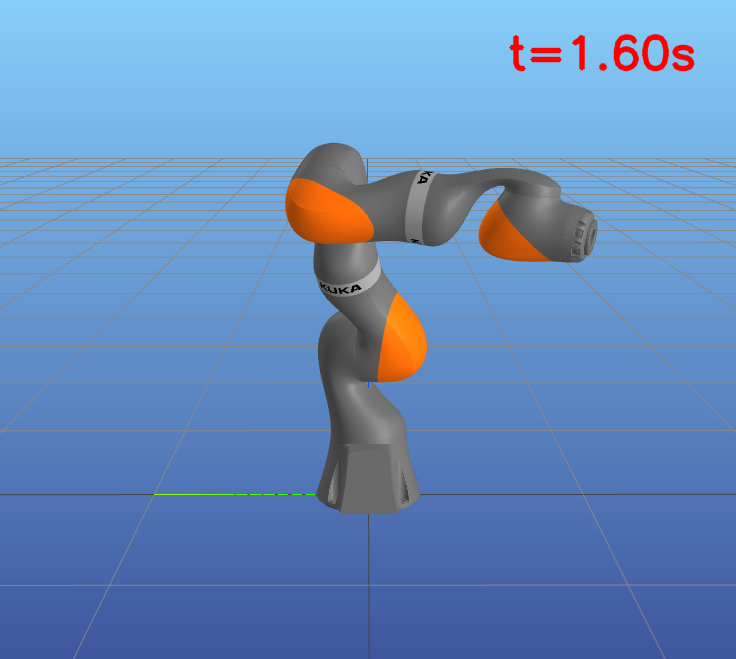}
    \includegraphics[width=0.15\textwidth]{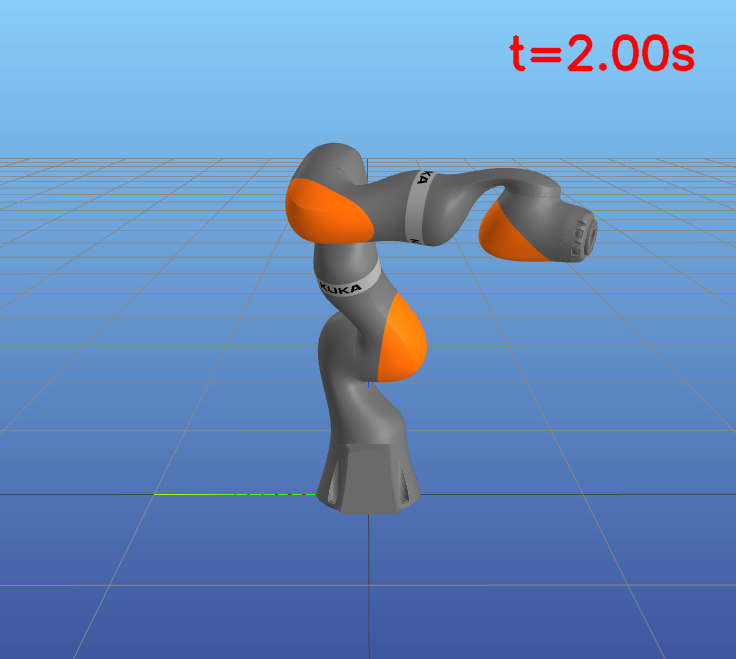}\\
    \includegraphics[width=0.15\textwidth]{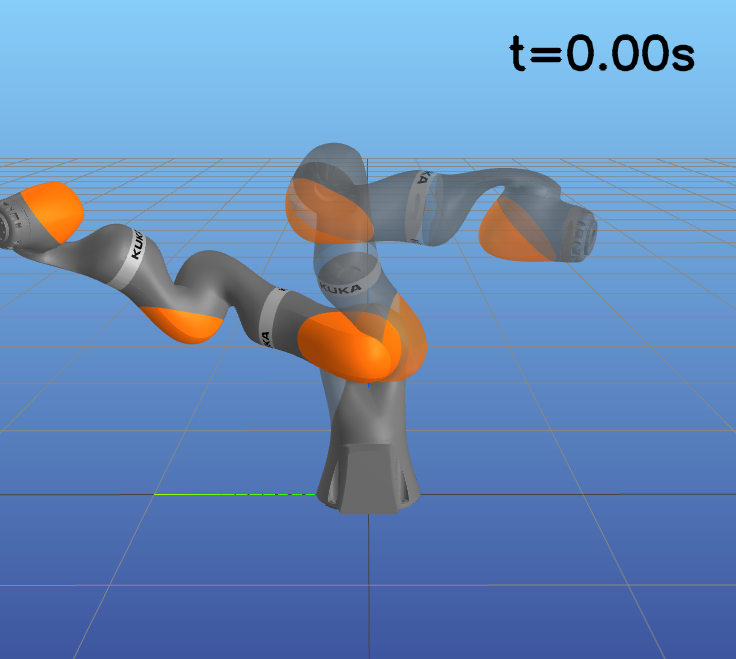}
    \includegraphics[width=0.15\textwidth]{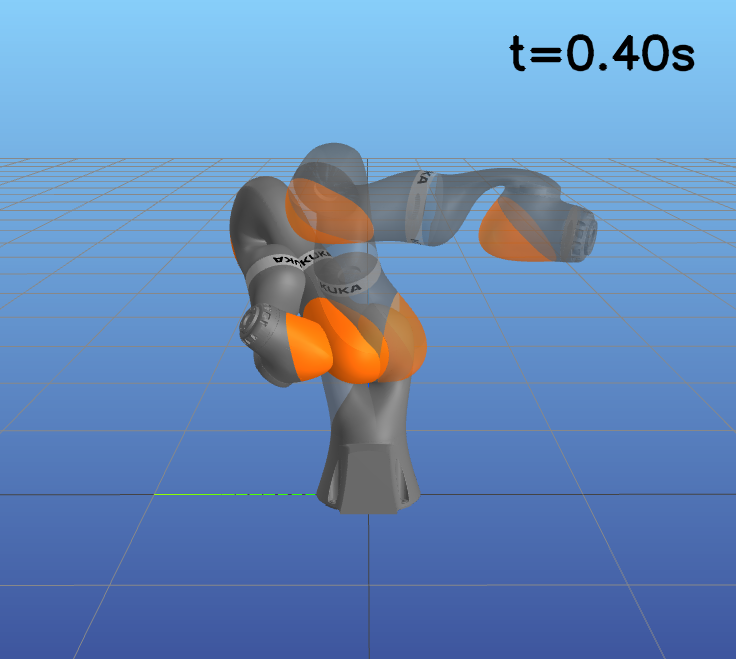}
    \includegraphics[width=0.15\textwidth]{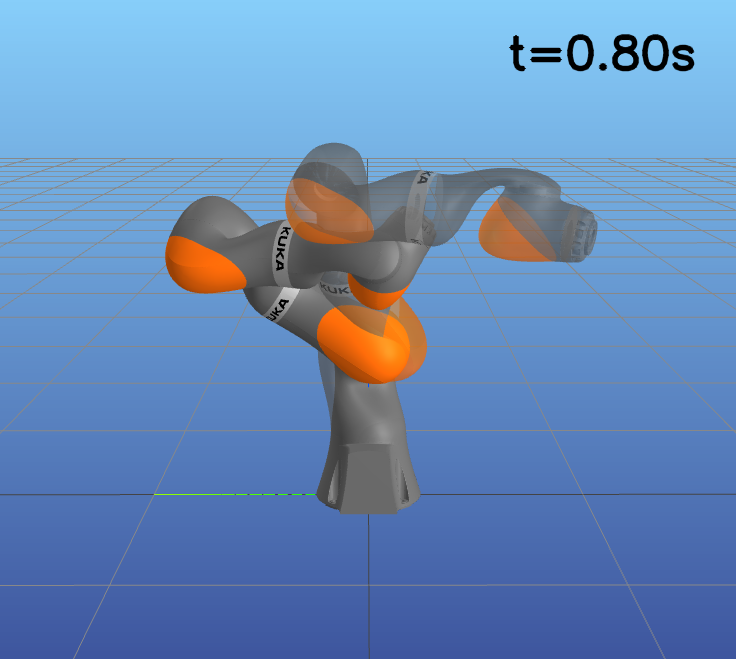}
    \includegraphics[width=0.15\textwidth]{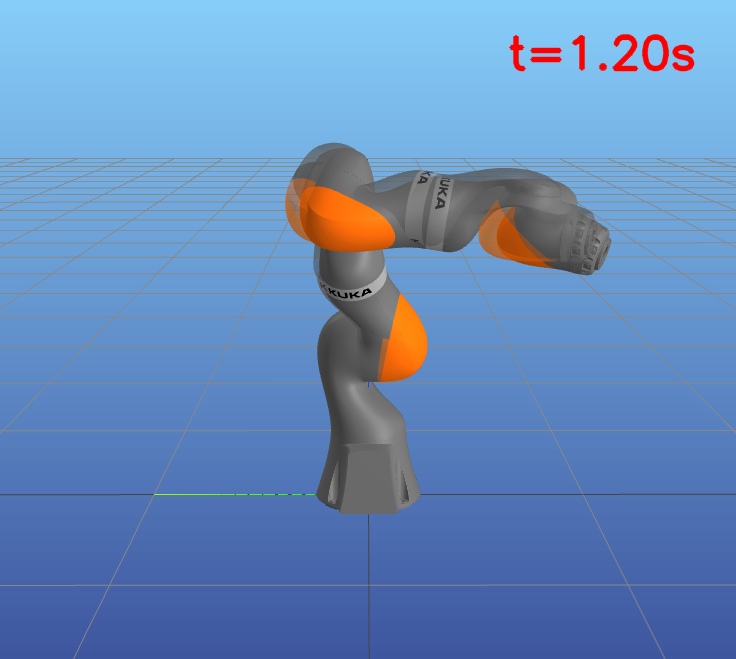}
    \includegraphics[width=0.15\textwidth]{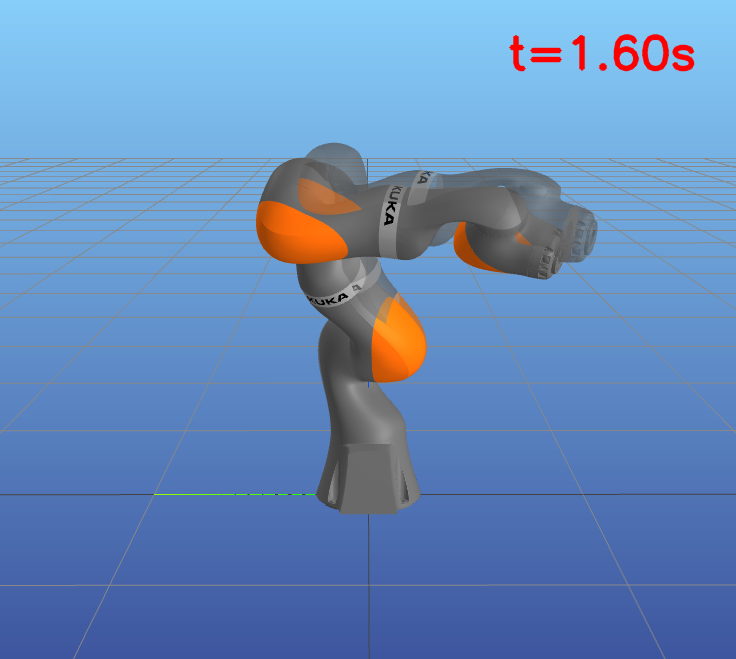}
    \includegraphics[width=0.15\textwidth]{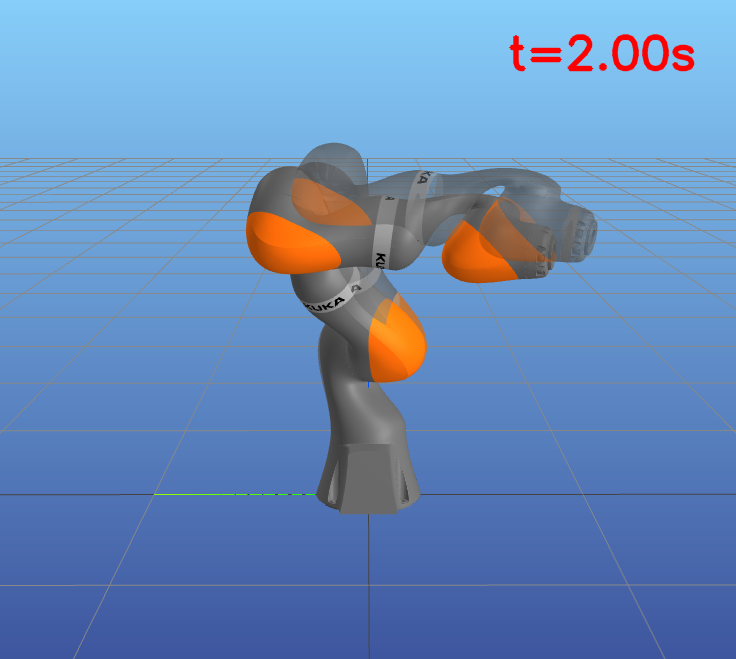}\\
    \includegraphics[width=0.15\textwidth]{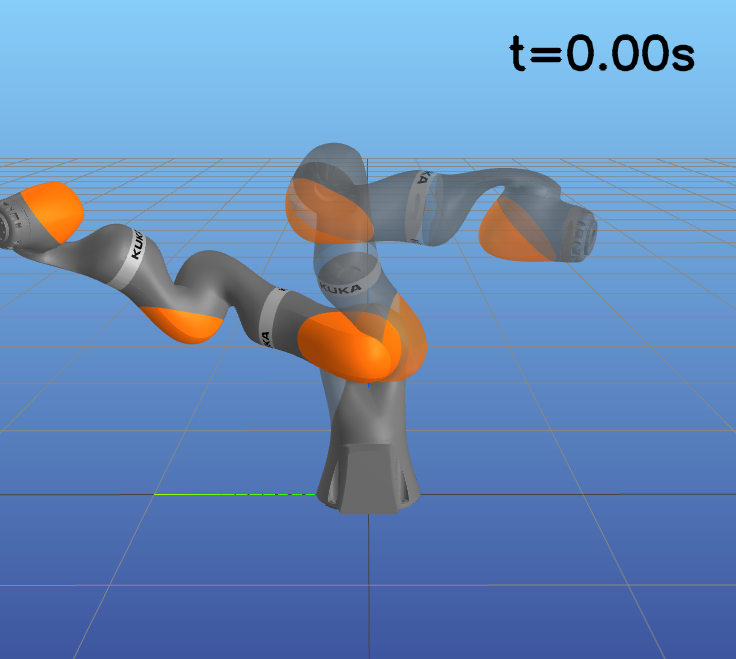}
    \includegraphics[width=0.15\textwidth]{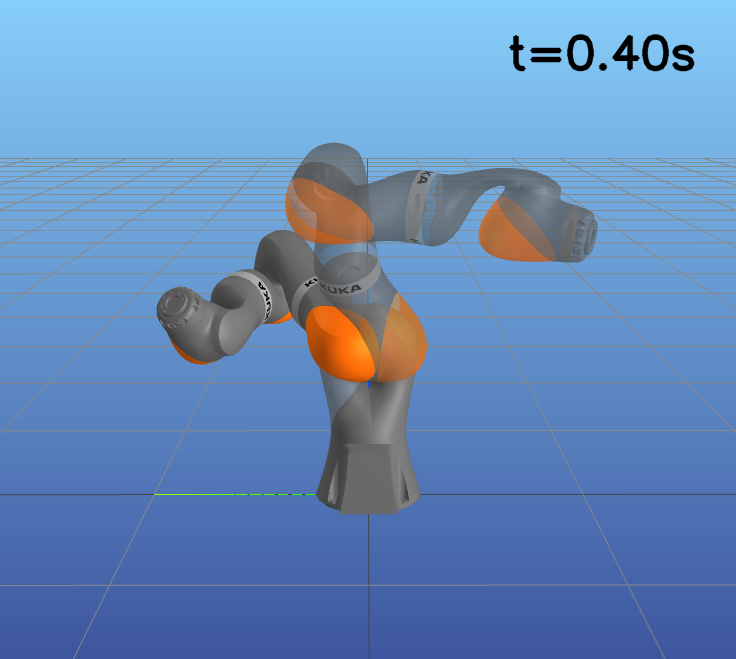}
    \includegraphics[width=0.15\textwidth]{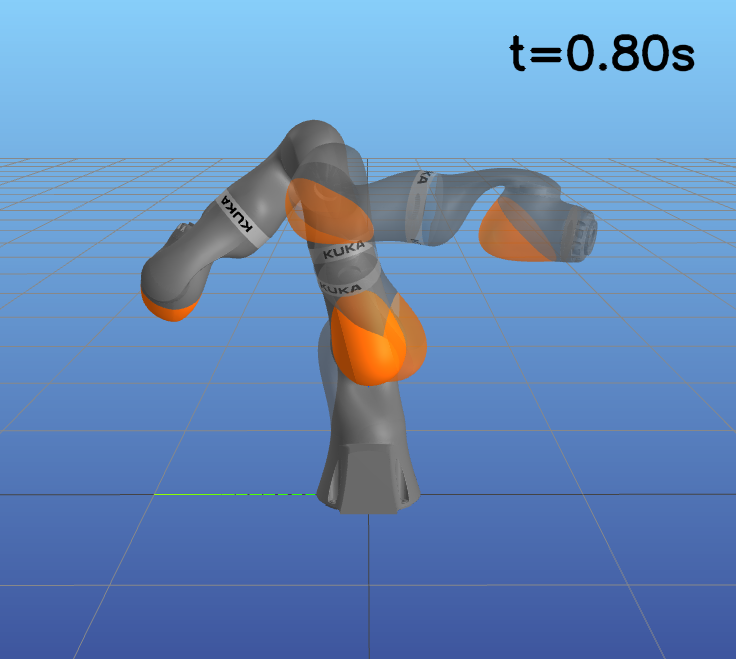}
    \includegraphics[width=0.15\textwidth]{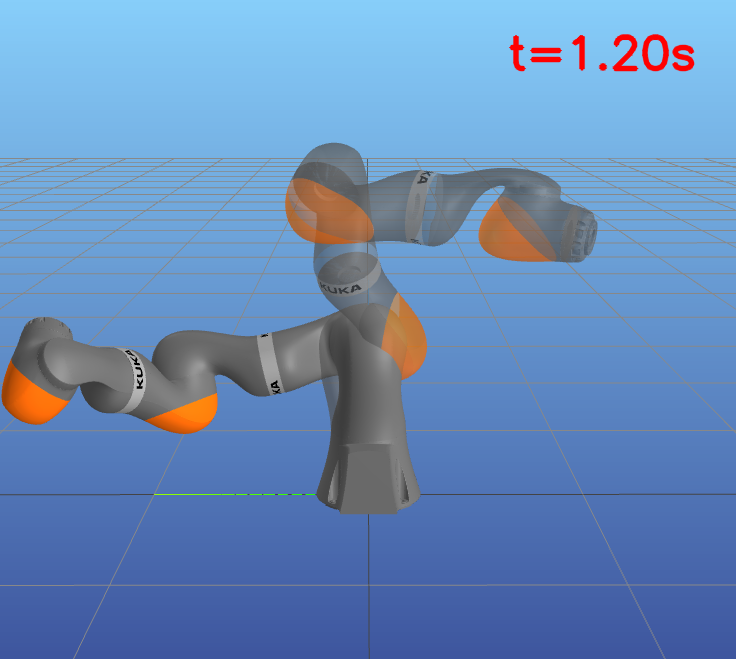}
    \includegraphics[width=0.15\textwidth]{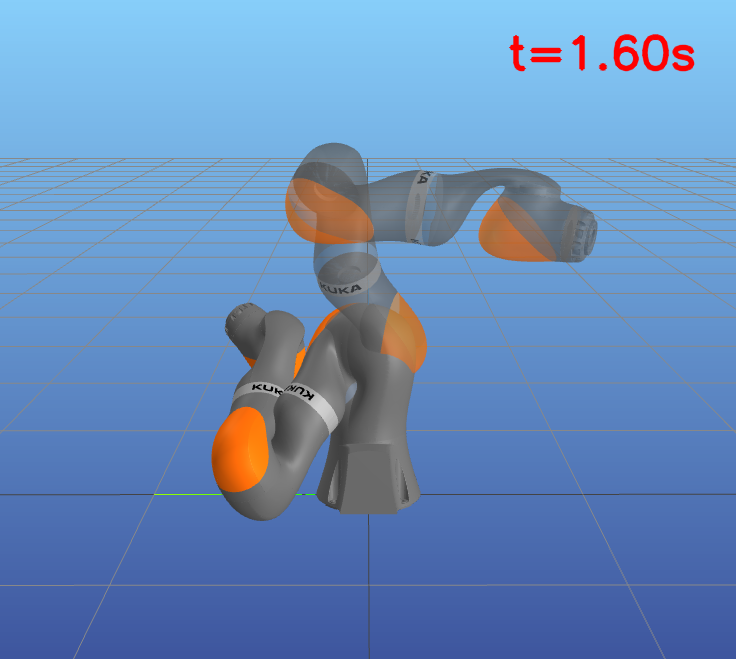}
    \includegraphics[width=0.15\textwidth]{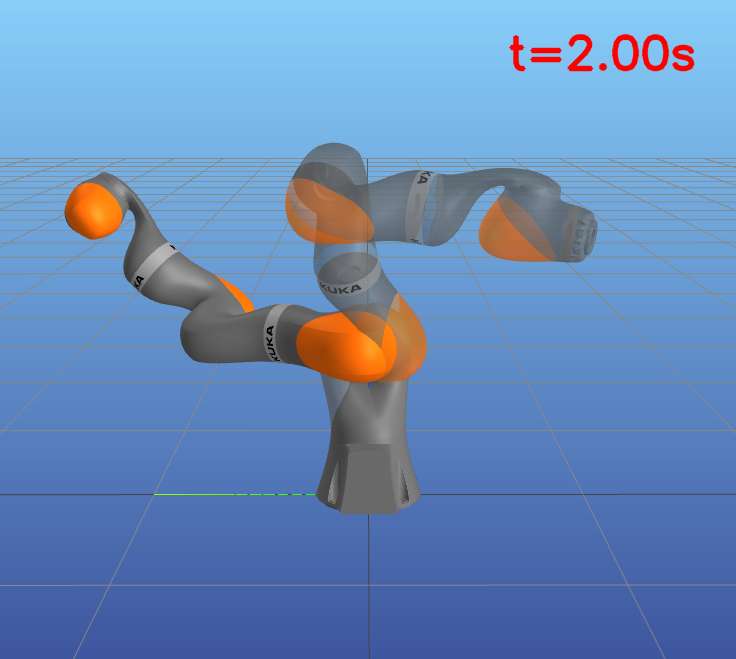}
    \caption{
    Simulation results of different controllers of different enhanced sampling strategies and different architectures.
    Top: IVP-ART with the modified QRnet ($\tau=1.0$) reaches the terminal state and maintains it;
    Middle: IVP-ART with MLP with ($\tau=1.0$) fails to reach the terminal state within an acceptable margin and exhibits instability around this state;
    Bottom: \DAgger\ with the modified QRnet does not approach the state at all in this particular failure case.
    The top row reaches the terminal state at the time of $1.1635$ seconds, after which, the timestamp color of all rows is changed to red.
    \label{fig:qr_demo_2}
    }
\end{figure*}

\begin{figure*}[h]
    \centering
    \includegraphics[width=0.15\textwidth]{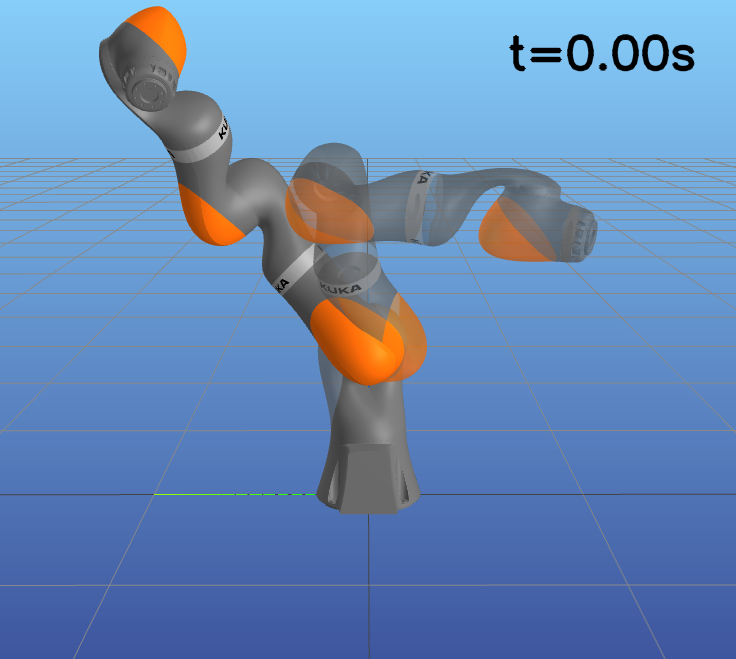}
    \includegraphics[width=0.15\textwidth]{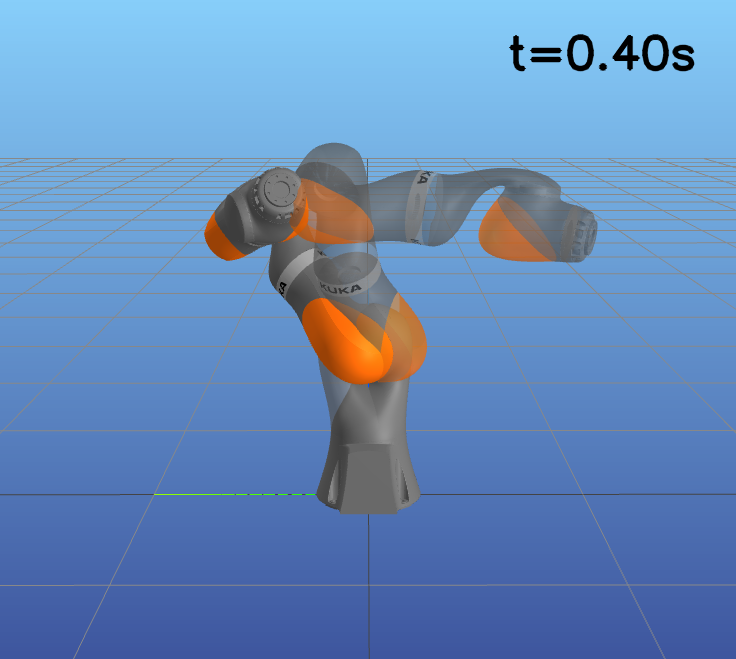}
    \includegraphics[width=0.15\textwidth]{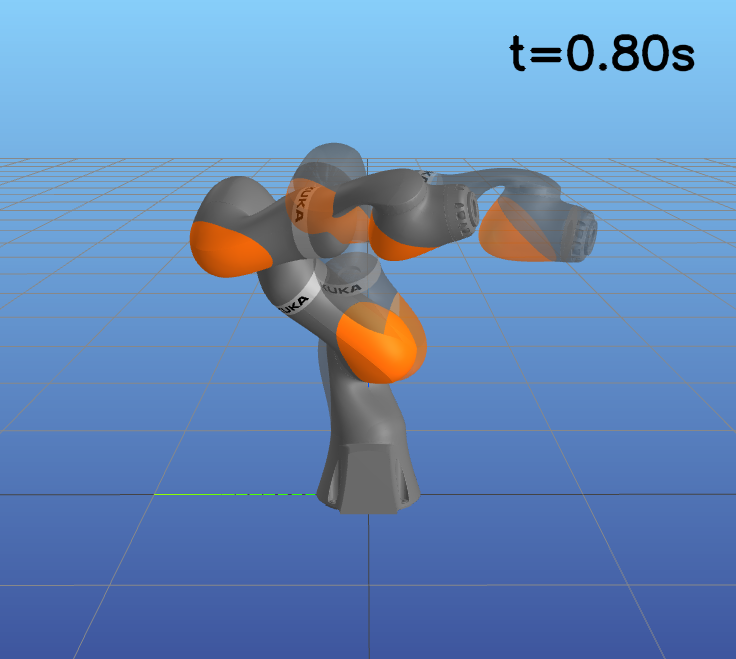}
    \includegraphics[width=0.15\textwidth]{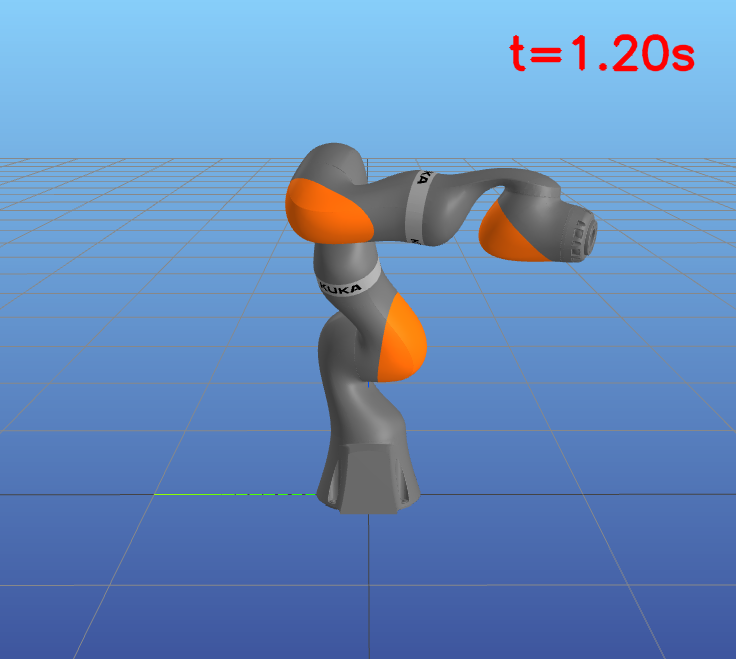}
    \includegraphics[width=0.15\textwidth]{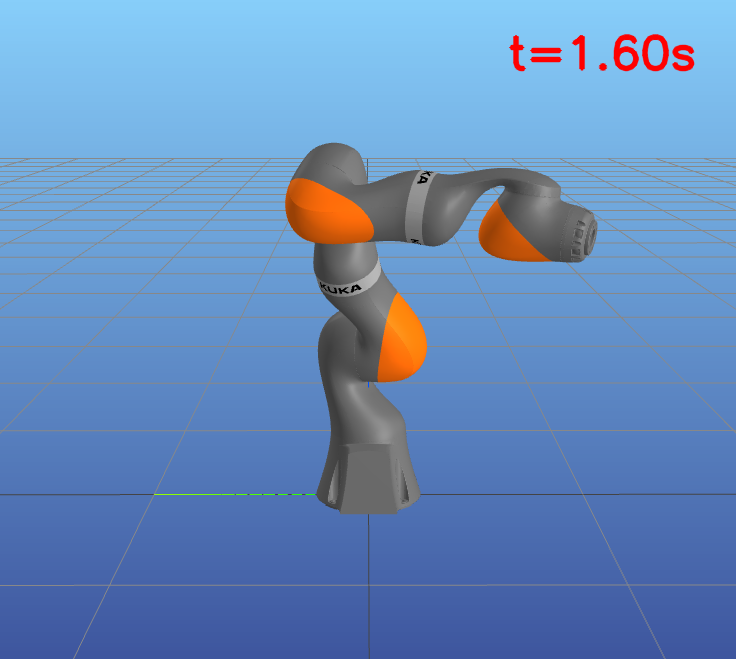}
    \includegraphics[width=0.15\textwidth]{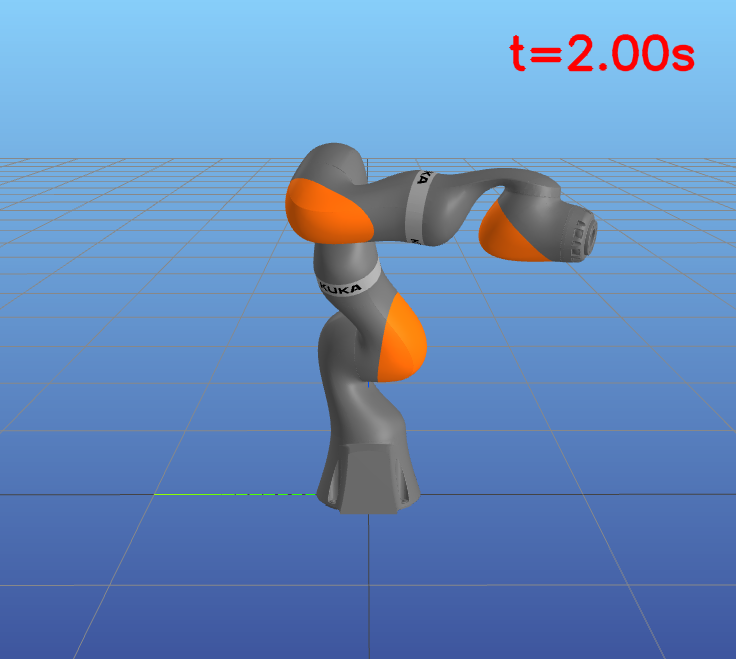}\\
    \includegraphics[width=0.15\textwidth]{fig/demos/8/qr/0.00s.png}
    \includegraphics[width=0.15\textwidth]{fig/demos/8/qr/0.40s.png}
    \includegraphics[width=0.15\textwidth]{fig/demos/8/qr/0.80s.png}
    \includegraphics[width=0.15\textwidth]{fig/demos/8/qr/1.20s.png}
    \includegraphics[width=0.15\textwidth]{fig/demos/8/qr/1.60s.png}
    \includegraphics[width=0.15\textwidth]{fig/demos/8/qr/2.00s.png}\\
    \includegraphics[width=0.15\textwidth]{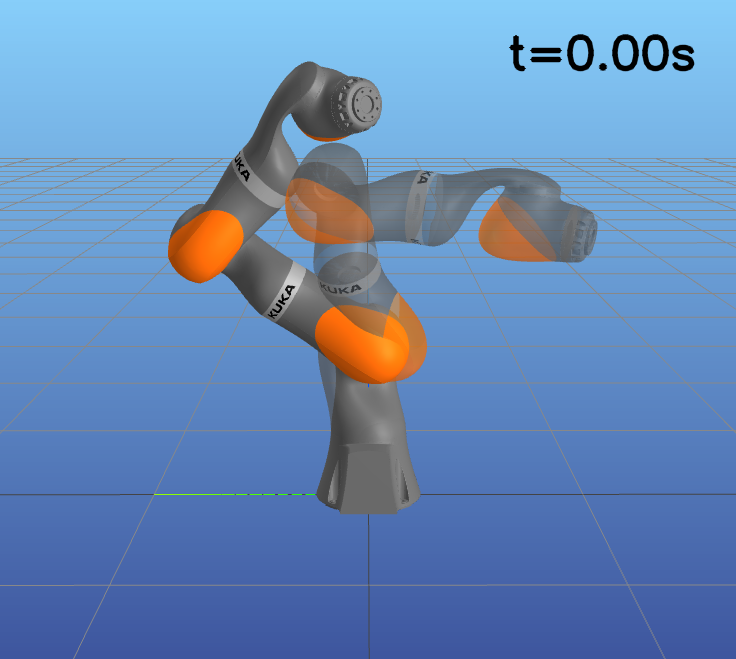}
    \includegraphics[width=0.15\textwidth]{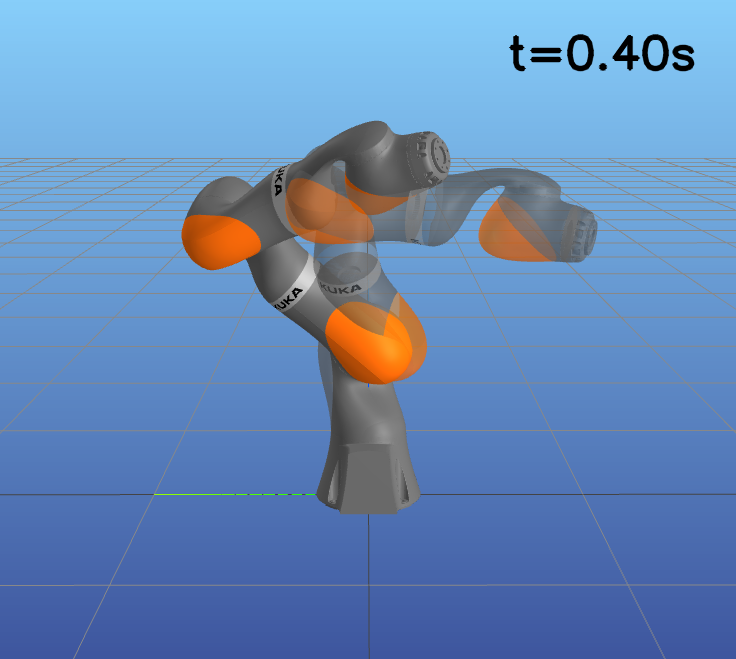}
    \includegraphics[width=0.15\textwidth]{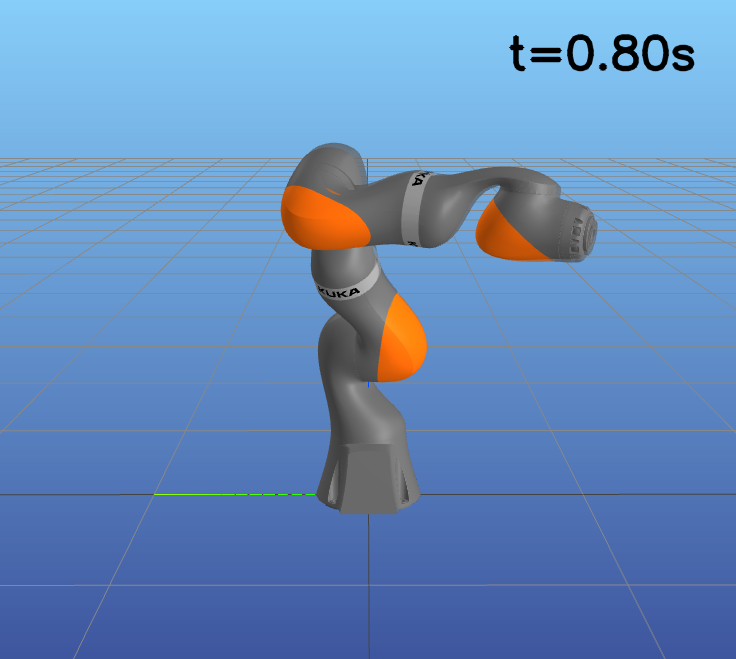}
    \includegraphics[width=0.15\textwidth]{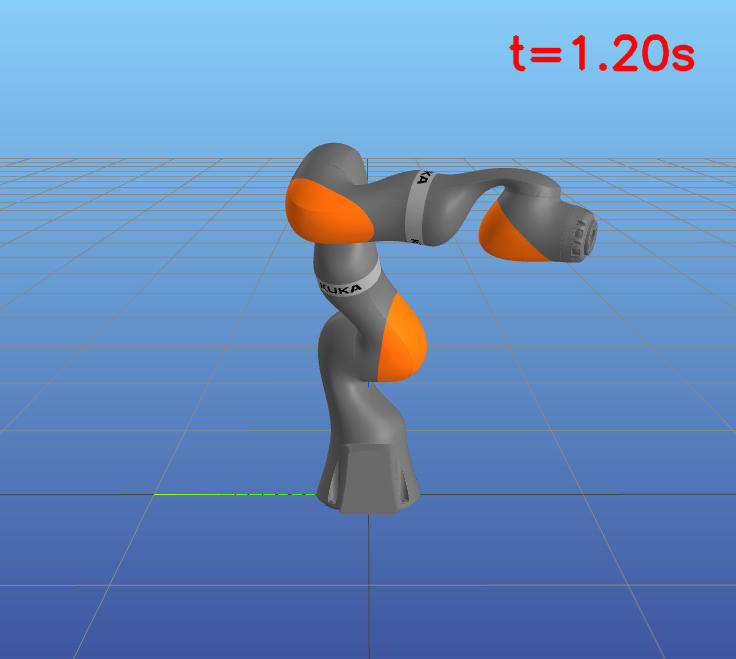}
    \includegraphics[width=0.15\textwidth]{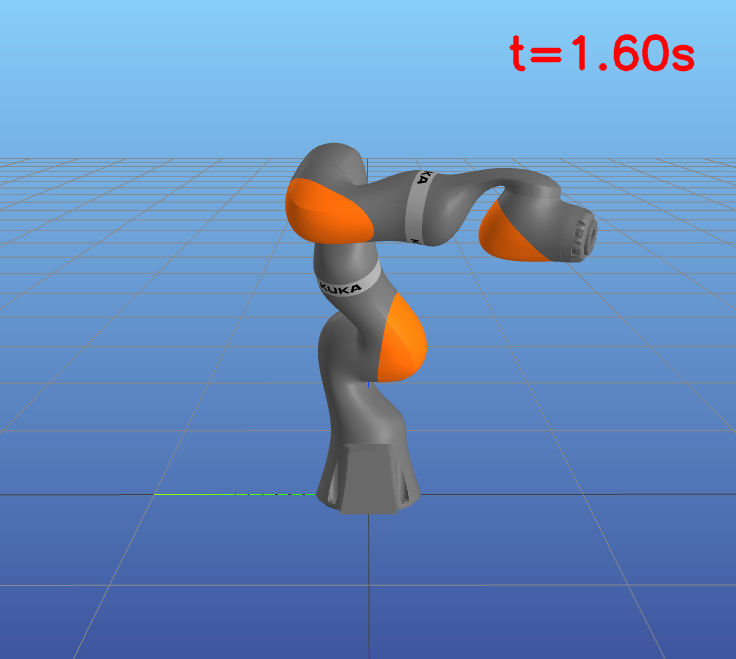}
    \includegraphics[width=0.15\textwidth]{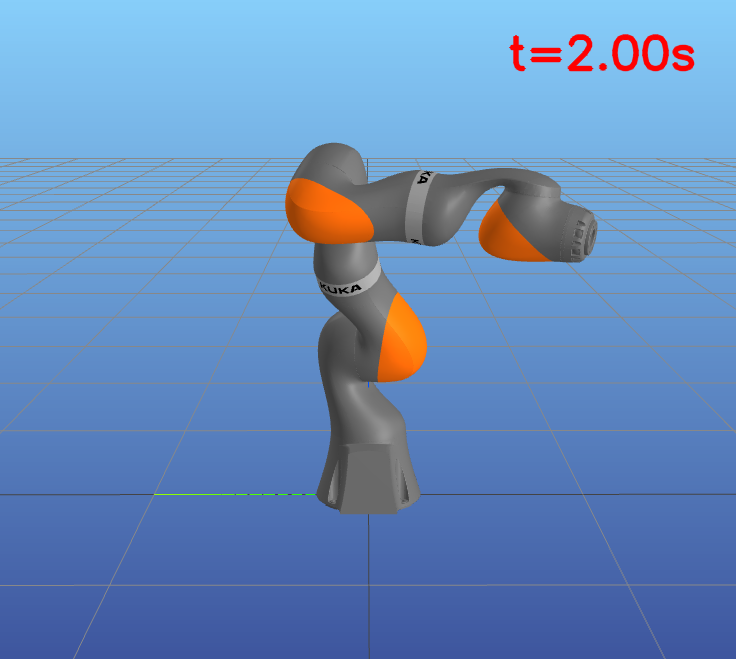}
    \caption{
    Simulation results of the controller of the IVP-ART algorithm, with $\tau=1.0$, applied to three distinct initial states, achieving terminal states at times $1.1815$, $1.1635$, and $0.8935$, respectively. Upon reaching the terminal state, the timestamp color is changed to red.
    }
    \label{fig:qr_demo_3}
\end{figure*}
   
\section{Demonstration}

This section presents visualizations of simulation trajectories at different time steps in Figure~\ref{fig:qr_demo_2} and Figure~\ref{fig:qr_demo_3}.

\end{document}